\PassOptionsToPackage{hypertexnames=false,breaklinks=true,bookmarks=true,bookmarksopen=false,draft=false}{hyperref}
\documentclass[final,hidelinks,onefignum,onetabnum]{siamart251216}

\usepackage{amsfonts,amssymb}
\usepackage{booktabs}
\usepackage{placeins}
\usepackage{enumitem}
\usepackage{tikz}
\hypersetup{hidelinks,colorlinks=true,urlcolor=blue,citecolor=blue,linkcolor=blue}

\newsiamremark{remark}{Remark}
\newsiamremark{example}{Example}
\crefname{example}{Example}{Examples}
\crefname{section}{Section}{Sections}
\crefname{subsection}{Subsection}{Subsections}
\newsiamthm{observation}{Observation}
\crefname{observation}{Observation}{Observations}

\makeatletter
\renewcommand*{\top}{%
  {\mathpalette\@transpose{}}%
}
\newcommand*{\@transpose}[2]{%
  \raisebox{\depth}{$\m@th#1\scriptscriptstyle\mathsf{T}$}%
}
\makeatother

\newcommand{\Po}{\mathcal{P}} 
\newcommand{\OP}{\mathcal{O}(P)} 
\newcommand{\D}{\mathcal{D}}
\DeclareMathOperator{\conv}{conv}
\DeclareMathOperator{\rank}{rank}
\DeclareMathOperator{\height}{height}
\DeclareMathOperator{\width}{width}
\DeclareMathOperator{\dir}{dir}
\DeclareMathOperator{\Span}{span}
\DeclareMathOperator*{\argmax}{argmax}

\newcommand{\SF}[1]{\cite[#1]{QL_staircase}}
\newcommand{\SFexample}{Example~3.2(ii)}%   alpha=(1,1,-1) for sums of two staircase simplices
\newcommand{\SFzeropadding}{Lemma~3.5}%    zero padding preserves the least number of summands
\newcommand{\SFsize}{Theorem~4.7}%         size of facet normals; (b) three summands
\newcommand{\SFcoeff}{Proposition~4.3}%    max entry over all sums = coeff(d) of 0/1-polytopes

\headers{On the facets of the Cayley polytope of Stanley's triangulation}{Y. Qu and J. Lee}

\title{On the facets of the Cayley polytope \\ of Stanley's triangulation\thanks{\today.
\funding{Supported in part by the Office of Naval Research grant N00014-24-1-2694.}
}}

\author{Yushan Qu\thanks{University of Michigan, Ann Arbor, MI
  (\email{yushanqu@umich.edu}, \email{jonxlee@umich.edu}).}
\and Jon Lee\footnotemark[2]}

\ifpdf
\hypersetup{
  pdftitle={On the facets of the Cayley polytope of Stanley's triangulation},
  pdfauthor={Y. Qu and J. Lee}
}
\fi

\begin{document}

\maketitle

% REQUIRED
\begin{abstract}
Disjunctions over the cells of a triangulation arise in formulations of piecewise-linear functions. For a poset $P$ on $[d]$, we study the convex hull $\D$ of the simplices of Stanley's triangulation of the order polytope of $P$, each lifted by its own indicator vector. 
Via the Cayley trick, for any disjunction of full-dimensional polytopes, the facets of the convex hull other than those described by the bounds on the indicator variables correspond to the facet normals of the Minkowski sum of the polytopes, and each is described by the full optimal big-M lifting of the supporting inequality with that normal of any one polytope.
For Stanley's triangulation, with $M(P)$ the Minkowski sum of its simplices, lifting only the facet-describing inequalities of the simplices need not suffice, unlike some other cases; we show that it suffices if and only if $M(P)$ is alcoved, if and only if deleting one element of $P$ leaves a chain. Otherwise, the coefficients of the additional facet-describing inequalities are governed by the least number $\ell(P)$ of elements whose deletion leaves a chain: for fixed $\ell(P)$, the largest entry of a primitive facet normal grows like $d^{\ell(P)-1}$, and it can be exponential in $d$ if $\ell(P)$ grows with $d$. For $d \geq 3$ and $\ell(P)\leq2$ the sharp bound is $d-2$.
\end{abstract}

% REQUIRED

% REQUIRED

\section{Introduction}\label{sec:intro}

For full-dimensional polytopes $\Po_0\,,\ldots,\Po_n\subset\mathbb{R}^d$ and indicator vectors $z\in\mathbb{R}^n$, where $z=\mathbf{e}_k$ (a standard unit vector) selects $\Po_k$ for $k\in N:=\{1,\ldots,n\}$ and $z=\mathbf{e}_0:=\mathbf{0}\in\mathbb{R}^n$ selects $\Po_0$\,, we study
\[
\D:=\conv\bigl(\textstyle\bigcup_{k\in N_0}(\Po_k\times\{\mathbf{e}_k\})\bigr)\subset\mathbb{R}^{d+n},\qquad N_0:=\{0,1,\ldots,n\}.
\]
The full-dimensional polytope $\D$ is the \emph{Cayley embedding}, or \emph{Cayley polytope}, of $\Po_0\,,\ldots,\Po_n$ (see \cite{HuberRambauSantos2000}, for example), a key tool for studying the Minkowski sum $\Po_0+\cdots+\Po_n$\,.

In a closely related stream of research motivated by devising models and algorithms for some real-world problems, \emph{disjunctive programming} studies optimization over a union of polyhedra (see \cite{BalasDP,balas2018disjunctive}). The convex hull of a union of polytopes has a well-known extended formulation (see \cite{BALAS1988279}; also \cite{Jeroslow,Blair}), but it introduces a copy of the variables for every polytope. Big-M formulations, which model the disjunction with indicator variables, are much more compact, but their continuous relaxations (with $z_k\in[0,1]$ in place of $z_k\in\{0,1\}$) can be quite weak; see \cite{Vielma2015MixedIL,TRESPALACIOS201598,Misener_BigM,GurobiM}. Disjunctive cuts and related strengthening techniques are used in many settings; see \cite{belotti2011disjunctive,Lodi-et-al:2020,kilincc2015two,Misener_Disj}.
Also motivated by disjunctive programming, \cite{Vielma2019Cayley} used the Cayley embedding, together with a non-polyhedral generalization of it, to obtain strong formulations for unions of convex sets without continuous auxiliary variables.

The facets of $\D$ that involve only $z$ are those of the simplex $\conv\{\mathbf{e}_0\,,\ldots,\mathbf{e}_n\}$, described by $z_k\geq0$ for $k\in N$ and $\sum_{k\in N}z_k\leq1$; the other facets are \emph{vertical}. A natural way to obtain a valid inequality for $\D$ is to start from a valid inequality $\alpha^\top x\leq\beta$ for one $\Po_k$ and to lift every indicator variable into it with the best coefficient that keeps it valid. In \cite{QL_ISCO2024,QU202531} we called the result the \emph{full optimal big-M lifting} of $\alpha^\top x\leq\beta$ (recalled in \cref{subsec:lifted}) and asked when the full optimal big-M liftings of the facet-describing inequalities of the $\Po_k$ describe all vertical facets of $\D$. The answer is positive for axis-aligned hyper-rectangles, for generalized cross polytopes, and whenever $d\in\{1,2\}$; see \cite{QU202531,qu2026disjunctionconvexhullsgeneralized}.

Here the polytopes are the cells of a triangulation. This setting arises when nonlinear functions are approximated by piecewise-linear ones, a standard way of bringing nonlinear functions within reach of mixed-integer linear optimization (see \cite{LeeWilson2001,Vielma2015MixedIL}): one triangulates the domain of a continuous function $f$, evaluates $f$ at the vertices, and interpolates linearly on each simplex. The interpolation is correct only if $x$ is a convex combination of the vertices of a \emph{single} simplex, which \cite{LeeWilson2001} enforces with one binary variable per simplex, these variables summing to one. Eliminating the variable of $\Po_0$ leaves exactly the indicator vectors $z$ above, and the convex hull of the feasible set in the $(x,z)$-space is $\D$.

We take the most basic family of triangulations with a combinatorial structure. For a poset $P$ on $[d]$, the order polytope $\OP\subset[0,1]^d$ (see \cite{Stanley1986}) is the union of the staircase simplices $\Delta_\sigma$ (\cref{subsec:order}) indexed by the linear extensions $\sigma$ of $P$, which form Stanley's triangulation of $\OP$. There can be as many as $d!$ linear extensions, so the extended formulation of $\D$ is large, and a description in the $(x,z)$-space is desirable. If $P$ is an antichain, Stanley's triangulation is the staircase triangulation of $[0,1]^d$, known in piecewise-linear approximation 
and fixed-point computation as the Freudenthal triangulation (see \cite{freudenthal1942simplizialzerlegungen}), the Kuhn subdivision (see \cite{Kuhn1960}), and Todd's triangulation $K_1$ (see \cite{todd1976computation}). Already for $d=3$, $\D$ then has vertical facets that are not described by the full optimal big-M lifting of a facet-describing inequality of any of the six simplices (\cref{ex:antichain}). We determine exactly when such facets appear and what they look like.

\smallskip

\noindent\emph{Contributions.}
By \cref{thm:verticalarelifting}, lifting facet-describing inequalities suffices if and only if $M(P)$ is alcoved, if and only if the simplices of Stanley's triangulation are pairwise compatible in the sense of \cite{early2025when}, if and only if the linear extension graph of $P$ is acyclic, if and only if deleting one element of $P$ leaves a chain; equivalently, $P$ has no induced subposet isomorphic to $A_3$\,, $2+2$, $\mathsf N$, or $K_{2,2}$ (\cref{cor:forbidden}). In particular, this happens only when Stanley's triangulation has at most $d$ simplices (\cref{rem:height2}).

This result rests on the correspondence between a Cayley embedding and the Minkowski sum of its summands (the Cayley trick). For any disjunction of full-dimensional polytopes, it gives a bijection between the vertical facets of the Cayley embedding and the facet-normal rays of the Minkowski sum, and we observe that the facet corresponding to $\alpha$ is described by the full optimal big-M lifting of the supporting inequality with normal $\alpha$ of \emph{any} one polytope (\cref{prop:cayley}). So every vertical facet is described by a lifted inequality; what can fail is that the inequality being lifted is facet-describing for the polytope: it may describe a lower-dimensional face. Hence the question is which facet normals the Minkowski sum has, and for Stanley's triangulation, where it is $M(P)$ (\cref{cor:bijection}), we answer it.

When lifting facet-describing inequalities does not suffice, the coefficients of the inequalities describing the additional facets are governed by the least number $\ell(P)$ of elements whose deletion leaves a chain. If $\ell(P)=0$, then $P$ is a chain and every primitive
facet normal of $M(P)$ has infinity norm one.
If $\ell(P)\ge 1$, primitive facet normals of $M(P)$ have entries at most $\delta_{\ell(P)}\binom{d-1}{\ell(P)-1}\leq\ell(P)\,(d-1)^{\ell(P)-1}$, where $\delta_\ell$ is the largest absolute determinant of an $\ell\times\ell$ zero-one matrix (\cref{thm:posetheight}), and for fixed $\ell\geq2$ this order is sharp: the largest entry $H_\ell(d)$ over posets with $\ell(P)=\ell$ is $\Theta(d^{\ell-1})$ (\cref{thm:product,cor:growth}). With $\ell(P)=\lfloor(d-1)/3\rfloor+1$ the entries can be of order $3^{d/3}$ (\cref{cor:exp}). For $\ell(P)=3$ the largest entry lies between $\lfloor(d-2)^2/3\rfloor-1$ and $\lfloor(d-1)^2/2\rfloor$ (\cref{cor:grid,thm:upper3}), and for $\ell(P)=4$ we have $1/27\leq\liminf_d H_4(d)/d^3\leq\limsup_d H_4(d)/d^3\leq1/2$ (\cref{cor:constants4}); the exact values and the sharp constants remain open (\cref{rem:open}). For $d \geq 3$ and $\ell(P)\leq2$ the sharp bound is $d-2$ (\cref{thm:sharpchaintwoelements}), attained by a chain with two isolated elements, for which we determine all vertical facets of $\D$ (\cref{cor:twoisolatescount}).

\smallskip
\noindent\emph{Outline.}
\Cref{sec:prelim} collects preliminaries on order polytopes, alcoved polytopes, restrictions of linear extensions, and support functions. \Cref{sec:vertical} formally defines full optimal big-M lifting, shows for arbitrary full-dimensional polytopes that the vertical facets of the Cayley embedding correspond to the facet normals of the Minkowski sum, specializes this to Stanley's triangulation, and proves \cref{thm:verticalarelifting}. \Cref{sec:nearchain} bounds the entries of the facet normals of $M(P)$ by a polynomial of degree $\ell(P)-1$ in $d$ and shows that this order is attained; it also classifies the facet normals for a chain with two isolated elements, shows that the entries can be exponential in $d$ when $\ell(P)$ grows with $d$, and sharpens the constants for $\ell(P)=3$ and $\ell(P)=4$. Minkowski sums of arbitrary staircase simplices, not tied to a poset, are the subject of the companion paper \cite{QL_staircase}, which studies how the number of summands, rather than the order relations, constrains their facet normals.

\smallskip
\noindent\emph{Notation.}
For a positive integer $d$, let $[d]:=\{1,2,\ldots,d\}$ and $[d_0]:=\{0,1,\ldots,d\}$.
We write $\mathbf{e}$ for the all-ones vector,  
and $\mathbf{e}_B:=\sum_{i\in B}\mathbf{e}_i$ for the indicator vector of a set $B$; the dimension will be clear from the context. We write $\lVert\alpha\rVert_\infty:=\max_i|\alpha_i|$. The symbol $\subset$ denotes inclusion, not necessarily proper (so $A\subset A$).

%%%%%%%%%%%%%%%%%%%%%%%%%%%

\section{Preliminaries}\label{sec:prelim}

\subsection{Order polytopes and linear extensions}\label{subsec:order}
Let $P=(V,\prec_P)$ be a poset on $V:=[d]$, where $\prec_P$ is a strict partial order; we read $i\prec_P j$ as ``$i$ lies below $j$''. Elements $i\neq j$ are \emph{comparable} if $i\prec_P j$ or $j\prec_P i$, and \emph{incomparable}, written $i\parallel_P j$, otherwise. A \emph{chain} is a set of pairwise comparable elements, labeled $c_1\prec_P\cdots\prec_P c_m$ (the empty set is a chain), and an \emph{antichain} is a set of pairwise incomparable elements; the \emph{height} $\height(P)$ and \emph{width} $\width(P)$ are the largest sizes of a chain and of an antichain. For $\varnothing\neq S\subset V$, the \emph{induced subposet} $P|_S$ is the restriction of $\prec_P$ to $S$, and for a sequence $\sigma$ of the elements of $V$, $\sigma|_S$ is the subsequence of elements of $S$. All notions below are defined in the same way for posets on any finite ground set $S$, with $\mathbb{R}^S$ in place of $\mathbb{R}^d$.

The \emph{order polytope} of $P$ is $\OP:=\{x\in[0,1]^d: x_i\leq x_j \text{ for } i\prec_P j\}$; it is a simplex if $P$ is a chain and the cube $[0,1]^d$ if $P$ is an antichain. A permutation $\sigma=(\sigma_1\,,\ldots,\sigma_d)$ of $V$ is a \emph{linear extension} of $P$ if $i\prec_P j$ implies that $i$ precedes $j$ in $\sigma$. We write $\mathcal L(P)=\{\sigma^0,\ldots,\sigma^n\}$, so $|\mathcal L(P)|=n+1$. For a permutation $\sigma$, the \emph{staircase simplex} $\Delta_\sigma:=\{x\in[0,1]^d: x_{\sigma_1}\leq\cdots\leq x_{\sigma_d}\}$ has volume $1/d!$, and the $d!$ staircase simplices triangulate $[0,1]^d$. Those contained in $\OP$ are exactly the $\Delta_\sigma$ with $\sigma\in\mathcal L(P)$, and they form \emph{Stanley's triangulation} $\mathcal T(P):=\{\Delta_k: k\in N_0\}$ of $\OP$, where $\Delta_k:=\Delta_{\sigma^k}$; see \cite{Stanley1986}.

The vertices of $\Delta_k$ are $v^k_r:=\sum_{t=r+1}^d\mathbf{e}_{\sigma^k_t}$ for $r\in[d_0]$, so $v^k_0=\mathbf{e}$ and $v^k_d=\mathbf{e}_0$\,, and $\Delta_k=\{x: -x_{\sigma^k_1}\leq0,\ x_{\sigma^k_d}\leq1,\ x_{\sigma^k_r}-x_{\sigma^k_{r+1}}\leq0 \text{ for } r\in[d-1]\}$. The gap coordinates $y_0=x_{\sigma^k_1}$\,, $y_r=x_{\sigma^k_{r+1}}-x_{\sigma^k_r}$ ($r\in[d-1]$), and $y_d=1-x_{\sigma^k_d}$ map $\Delta_k$ affinely onto the standard simplex $\Delta^{\mathrm{std}}_d:=\{y\in\mathbb{R}^{d+1}_+:\sum_r y_r=1\}$, with inverse $x_{\sigma^k_i}=\sum_{r<i}y_r$; the vertex $v^k_r$ corresponds to $y_r=1$.

Let $G(P)$ be the graph on $\mathcal L(P)$ in which $\sigma\neq\tau$ are adjacent when $\dim(\Delta_\sigma\cap\Delta_\tau)=d-1$, that is, when the two simplices share a facet. By \cref{lem:adjacentswap} below, this happens exactly when $\tau$ is obtained from $\sigma$ by interchanging two incomparable elements in adjacent positions, so $G(P)$ is the \emph{linear extension graph} of $P$ (see \cite{Mas09}).

\subsection{Alcoved polytopes}\label{subsec:alcoved}
Following \cite{early2025when}, we use alcoved polytopes, written in our $x$-coordinates. A full-dimensional polytope $Q\subset\mathbb{R}^d$ is \emph{alcoved} if every facet normal of $Q$ is a positive multiple of one of the vectors $\mathbf{e}_i$\,, $-\mathbf{e}_i$\,, or $\mathbf{e}_i-\mathbf{e}_j$ ($i,j\in[d]$, $i\neq j$); equivalently, $Q$ can be described by inequalities of the forms $x_i\leq b_{i0}$\,, $-x_i\leq b_{0i}$ and $x_i-x_j\leq b_{ij}$ with $b_{ij}\in\mathbb{R}$.
By the description of $\Delta_k$ in \cref{subsec:order}, each simplex $\Delta_k$ is alcoved, and so is the order polytope $\OP$\,.
Two alcoved polytopes $Q_1$ and $Q_2$ are \emph{compatible} if their Minkowski sum $Q_1+Q_2$ is alcoved (see \cite{early2025when}).

\smallskip
\noindent\emph{Relation to \cite{early2025when}.}
Index the coordinates of $\mathbb{R}^{d+1}$ by $[d_0]$ and let $\mathcal{H}_{d+1} := \{w \in \mathbb{R}^{d+1}: \sum_{i \in [d_0]} w_i = 0\}$. In \cite{early2025when}, a polytope in $\mathcal{H}_{d+1}$ is alcoved if all its facet normals are parallel to the roots of type $A_d$\,, that is, to the vectors $\mathbf{e}_i-\mathbf{e}_j\in\mathbb{R}^{d+1}$, $i\neq j\in[d_0]$\,; the defining parameters are real. Define the linear isomorphism $\psi: \mathbb{R}^d \rightarrow \mathcal{H}_{d+1}$ by
\[
[\psi(x)]_0 := -\bar x\,, \qquad [\psi(x)]_i := x_i-\bar x\,, \qquad \bar x:=\textstyle\sum_{j \in [d]} x_j/(d+1)\,, \quad i \in [d],
\]
the orthogonal projection of $(0,x)$ onto $\mathcal{H}_{d+1}$\,. Because $[\psi(x)]_i-[\psi(x)]_0=x_i$ and $[\psi(x)]_i-[\psi(x)]_j=x_i-x_j$ for $i,j\in[d]$, a polytope $Q\subset\mathbb{R}^d$ is alcoved in the sense above if and only if $\psi(Q)$ is alcoved in the sense of \cite{early2025when}; and because $\psi$ is linear, $\psi(Q_1+Q_2)=\psi(Q_1)+\psi(Q_2)$, so compatibility is also preserved.

For $\sigma\in\mathcal{L}(P)$, the description of $\Delta_\sigma$ gives
\begin{equation}\label{eq:EKMdict}
\psi(\Delta_\sigma)=\{w\in\mathcal{H}_{d+1}:\ w_{\sigma_d}\geq w_{\sigma_{d-1}}\geq\cdots\geq w_{\sigma_1}\geq w_0\geq w_{\sigma_d}-1\}.
\end{equation}
This is the alcoved simplex $\Delta_{\mathbf S}$ of \cite[Definition 4.1]{early2025when} for the nondegenerate ordered set partition $\mathbf{S}_\sigma:=(\{\sigma_d\},\{\sigma_{d-1}\},\ldots,\{\sigma_1\},\{0\})$ of $[d_0]$ (after relabeling $[d_0]$ as $[d+1]$). Thus $\Delta_\sigma$ corresponds to the cyclic order $(\sigma_d\,,\ldots,\sigma_1\,,0)$ on $[d_0]$, and $\Delta_\sigma$ and $\Delta_\tau$ are compatible if and only if $\mathbf S_\sigma$ and $\mathbf S_\tau$ are compatible in the sense of \cite{early2025when}. The same holds, after relabeling, for linear extensions of a poset on any finite ground set $S$ in place of $[d]$, with the element $0$ adjoined; we use this for induced subposets $P|_S$\,.

\subsection{Restrictions and adjacent linear extensions}\label{subsec:posetlemmas}
\begin{lemma}\label{lem:restriction}
    Let $S \subset V$ be nonempty, and let $\pi_S: \mathbb{R}^V \rightarrow \mathbb{R}^S$, $\pi_S(x) = (x_s)_{s \in S}$\,, be the coordinate projection. Then:
    \begin{enumerate}
        \item The restriction $\sigma\mapsto\sigma|_S$ is a well-defined surjection $\mathcal{L}(P)\rightarrow\mathcal{L}(P|_S)$.
        \item For every $\sigma \in \mathcal{L}(P)$, $\pi_S(\Delta_{\sigma}) = \Delta_{{\sigma}|_S}$\,.
        \item If $\Delta_{\sigma} + \Delta_{\tau}$ is alcoved for all $\sigma, \tau \in \mathcal{L}(P)$, then $\Delta_{\rho} + \Delta_{\eta}$ is alcoved for all $\rho, \eta \in \mathcal{L}(P|_S)$.
    \end{enumerate}
\end{lemma}

\begin{proof}
    \emph{Item 1.} Let $\sigma \in \mathcal{L}(P)$. If $a,b \in S$ and $a \prec_{P|_S} b$, then $a \prec_{P} b$, so $a$ appears before $b$ in $\sigma$; deleting the elements of $V \setminus S$ does not change the relative order of $a$ and $b$. Hence $\sigma|_S \in \mathcal{L}(P|_S)$.

    For surjectivity, let $\rho=(\rho_1\,,\ldots,\rho_m) \in \mathcal{L}(P|_S)$. Let $G_\rho$ be the directed graph on $V$ with an arc $u \rightarrow v$ whenever $u \prec_P v$, and an arc $\rho_i \rightarrow \rho_j$ whenever $i < j$. We claim that $G_\rho$ is acyclic. A directed cycle cannot use only arcs of the first kind, because $P$ is a poset; so it contains an arc of the second kind, and hence elements of $S$. Split the cycle at its elements of $S$ into consecutive segments $s\to\cdots\to s'$ whose interior vertices lie outside $S$. A segment is either a single arc $\rho_i\to\rho_j$ with $i<j$, or a directed path of arcs of the first kind (arcs of the second kind join two elements of $S$), in which case $s\prec_P s'$ by transitivity, so $s\prec_{P|_S}s'$ and $s$ precedes $s'$ in $\rho$. In both cases the position in $\rho$ strictly increases along each segment, which is impossible around a cycle. Thus $G_\rho$ is acyclic, and any topological ordering $\sigma$ of $G_\rho$ is a linear extension of $P$ with $\sigma|_S=\rho$.

    \emph{Item 2.} If $x\in\Delta_\sigma$\,, its coordinates lie in $[0,1]$ and are weakly increasing along $\sigma$, hence along $\sigma|_S$\,, so $\pi_S(x)\in\Delta_{\sigma|_S}$\,. Conversely, let $y\in\Delta_{\sigma|_S}$ and define $x\in\mathbb{R}^V$ by $x_v:=\max\bigl(\{0\}\cup\{y_s:\ s\in S \text{ is } v \text{ or appears before } v \text{ in } \sigma\}\bigr)$. Then $x\in[0,1]^V$, the coordinates of $x$ are weakly increasing along $\sigma$, so $x\in\Delta_\sigma$\,, and $x_s=y_s$ for $s\in S$ because $y$ is weakly increasing along $\sigma|_S$\,. Hence $\pi_S(x)=y$.

    \emph{Item 3.} We first show that the coordinate projection of a full-dimensional alcoved polytope $Q\subset\mathbb{R}^V$ is alcoved. Introduce an auxiliary coordinate $x_0:=0$. Then $Q$ is described by inequalities $x_i-x_j\leq b_{ij}$ with $i\neq j\in V\cup\{0\}$. Eliminating one coordinate $x_v$\,, $v\in V$, by Fourier--Motzkin elimination combines each pair $x_i-x_v\leq b_{iv}$\,, $x_v-x_j\leq b_{vj}$ into $x_i-x_j\leq b_{iv}+b_{vj}$ (a constant inequality if $i=j$), so the projection is again described by difference inequalities. Iterating over $v\in V\setminus S$ shows that $\pi_S(Q)$ is described by inequalities of the forms $x_i-x_j\leq b$, $x_i\leq b$, $-x_i\leq b$ with $i,j\in S$; it is full dimensional because $Q$ is, so it is alcoved.

    Now let $\rho,\eta\in\mathcal L(P|_S)$. By item 1, there are $\sigma,\tau\in\mathcal L(P)$ with $\sigma|_S=\rho$ and $\tau|_S=\eta$. Because $\pi_S$ is linear, item 2 gives $\pi_S(\Delta_\sigma+\Delta_\tau)=\Delta_\rho+\Delta_\eta$\,. By hypothesis $\Delta_\sigma+\Delta_\tau$ is alcoved, so $\Delta_\rho+\Delta_\eta$ is alcoved.
\end{proof}

\begin{lemma}\label{lem:adjacentswap}
    Two distinct linear extensions $\sigma,\tau\in \mathcal{L}(P)$ are adjacent in $G(P)$ if and only if one is obtained from the other by interchanging two incomparable elements in adjacent positions.
\end{lemma}

\begin{proof}
    Let $\sigma\neq\tau$ in $\mathcal L(P)$, and let $\operatorname{Inv}(\sigma,\tau)$ be the set of pairs $\{a,b\}\subset V$ that $\sigma$ and $\tau$ place in opposite orders; it is nonempty because $\sigma\neq\tau$.

    We first claim that $\Delta_\sigma\cap\Delta_\tau=\{x\in\Delta_\sigma:\ x_a=x_b \text{ for all } \{a,b\}\in\operatorname{Inv}(\sigma,\tau)\}$. If $x\in\Delta_\sigma\cap\Delta_\tau$ and $\{a,b\}\in\operatorname{Inv}(\sigma,\tau)$, then one simplex gives $x_a\leq x_b$ and the other $x_b\leq x_a$\,, so $x_a=x_b$\,. Conversely, let $x\in\Delta_\sigma$ satisfy these equations, and let $u,v$ be consecutive in $\tau$, with $u$ first. If $\sigma$ also places $u$ before $v$, then $x_u\leq x_v$ because $x\in\Delta_\sigma$; otherwise $\{u,v\}\in\operatorname{Inv}(\sigma,\tau)$ and $x_u=x_v$\,. Together with $x\in[0,1]^d$, this gives $x\in\Delta_\tau$\,.

    For $x\in\Delta_\sigma$ the coordinates $x_{\sigma_1}\leq\cdots\leq x_{\sigma_d}$ are weakly increasing, so for $p<q$ we have $x_{\sigma_p}=x_{\sigma_q}$ if and only if $x_{\sigma_i}=x_{\sigma_{i+1}}$ for every $i$ with $p\leq i<q$. Let $J:=\bigl\{i\in[d-1]:\ p\leq i<q \text{ for some } \{\sigma_p\,,\sigma_q\}\in\operatorname{Inv}(\sigma,\tau) \text{ with } p<q\bigr\}$. By the claim, $\Delta_\sigma\cap\Delta_\tau=\Delta_\sigma\cap\{x:\ x_{\sigma_i}=x_{\sigma_{i+1}} \text{ for all } i\in J\}$, which is the face of the $d$-simplex $\Delta_\sigma$ obtained by making $|J|$ of its facet-describing inequalities tight; hence $\dim(\Delta_\sigma\cap\Delta_\tau)=d-|J|$. In particular, $\sigma$ and $\tau$ are adjacent if and only if $|J|=1$.

    Suppose that $J=\{i\}$. Every $\{\sigma_p\,,\sigma_q\}\in\operatorname{Inv}(\sigma,\tau)$ with $p<q$ satisfies $\{p,\ldots,q-1\}\subset J$, so $p=i$ and $q=i+1$. Thus the only pair that $\sigma$ and $\tau$ order differently is $\{\sigma_i\,,\sigma_{i+1}\}$, and because a permutation is determined by the relative orders of its pairs, $\tau$ is obtained from $\sigma$ by interchanging $\sigma_i$ and $\sigma_{i+1}$\,. These two elements are incomparable in $P$, because the linear extensions $\sigma$ and $\tau$ order them oppositely.

    Conversely, suppose that $\tau$ is obtained from $\sigma$ by interchanging two incomparable elements $\sigma_i$ and $\sigma_{i+1}$ in adjacent positions. The interchange changes only the relative order of this incomparable pair, so $\tau\in\mathcal L(P)$, $\operatorname{Inv}(\sigma,\tau)=\{\{\sigma_i\,,\sigma_{i+1}\}\}$, and $J=\{i\}$. Hence $\dim(\Delta_\sigma\cap\Delta_\tau)=d-1$, and $\sigma$ and $\tau$ are adjacent.
\end{proof}

\begin{lemma}\label{lem:posetchain}
    Let $P$ be a poset on a finite set $V$.
    \begin{enumerate}[label=(\alph*)]
        \item Suppose that $P$ has no three pairwise incomparable elements. If $a_1\parallel_P a_2$ and $a_3\parallel_P a_4$ with $a_1\,,\ldots,a_4$ distinct, then $\{a_1\,,\ldots,a_4\}=\{b_1\,,b_2\,,t_1\,,t_2\}$ with $b_1\parallel_P b_2$\,, $t_1\parallel_P t_2$\,, and $t_i\not\prec_P b_j$ for all $i,j$.
        \item Suppose that $P$ has no three pairwise incomparable elements. If any two incomparable pairs share an element, then $V\setminus\{y\}$ is a chain for some $y\in V$.
        \item If $V\setminus\{y\}$ is not a chain for any $y\in V$, then $P$ has three pairwise incomparable elements, or there are four distinct elements $b_1\,,b_2\,,t_1\,,t_2\in V$ with $b_1\parallel_P b_2$\,, $t_1\parallel_P t_2$\,, and $t_i\not\prec_P b_j$ for all $i,j$.
    \end{enumerate}
\end{lemma}

\begin{proof}
    (a) Let $S=\{a_1\,,a_2\,,a_3\,,a_4\}$. Every three-element subset of $S$ contains $\{a_1\,,a_2\}$ or $\{a_3\,,a_4\}$, so $P|_S$ contains no three-element chain, and by hypothesis no three-element antichain. Hence $P|_S$ has exactly two minimal elements: three minimal elements would form an antichain, and a unique minimal element would lie below the other three, which would then be pairwise incomparable (a relation between two of them would give a three-element chain). Let $b_1\,,b_2$ be the two minimal elements, which are incomparable, and $t_1\,,t_2$ the remaining elements. If, say, $t_1\prec_P t_2$\,, then a minimal element below $t_1$ together with $t_1\,,t_2$ would form a three-element chain; so $t_1\parallel_P t_2$\,. Finally, no $t_i$ lies below a $b_j$ because the $b_j$ are minimal in $P|_S$\,.

    (b) If there are no incomparable pairs, choose any $y\in V$; if there is exactly one, choose either of its elements. Otherwise, choose two distinct incomparable pairs; because they intersect, write them as $\{y,a\}$ and $\{y,b\}$. An incomparable pair not containing $y$ would have to meet both of these pairs, so it would be $\{a,b\}$; but then $y,a,b$ would be pairwise incomparable, which is excluded. Therefore every incomparable pair contains $y$, so any two distinct elements of $V\setminus\{y\}$ are comparable.

    (c) Suppose that $P$ has no three pairwise incomparable elements. Because $V\setminus\{y\}$ is not a chain for any $y$, (b) shows that some two incomparable pairs are disjoint, and (a) applied to them gives the four elements $b_1\,,b_2\,,t_1\,,t_2$\,.
\end{proof}

\subsection{Support functions and maximizing blocks}\label{subsec:support}
We use standard terminology for polytopes (see \cite{ziegler,Grunbaum2003}). An inequality $\alpha^\top x\leq\beta$ is \emph{valid} for a polytope $K$ if every point of $K$ satisfies it, and it \emph{describes} the face $K\cap\{x:\alpha^\top x=\beta\}$. A full-dimensional polytope has a minimal description by facet-describing inequalities, unique up to positive scaling.

For a nonempty polytope $K \subset \mathbb{R}^d$, the \emph{support function} is $h_K(\alpha) := \max_{x \in K} \alpha^\top x$, and the \emph{supporting inequality} $\alpha^\top x \leq h_K(\alpha)$ describes the nonempty face $F_K(\alpha) := \{x \in K: \alpha^\top x = h_K(\alpha)\}$; note $F_K(0)=K$.
If $K$ is full dimensional, a vector $\alpha\neq 0$ is a \emph{facet normal} of $K$ if $\dim F_K(\alpha)=d-1$; thus facet normals are outward normals. The \emph{facet-normal ray} of a facet normal $\alpha$ is $\{\lambda \alpha: \lambda > 0\}$, so two facet normals define the same facet-normal ray if and only if one is a positive multiple of the other.
The \emph{direction space} of a nonempty convex set $C$ is $\dir C:=\Span\{x-y: x,y\in C\}$, so $\dim C=\dim\dir C$, and $\dir(A+B)=\dir A+\dir B$ for nonempty convex sets $A,B$.

\smallskip
\noindent\emph{The Minkowski sum $M(P)$.}
For $k\in N_0$\,, let $h_k(\alpha) := h_{\Delta_k}(\alpha)$, and define $M(P) := \sum_{k \in N_0} \Delta_{k}$\,. Support functions and faces behave additively under Minkowski sums (see \cite{Schneider2014}): for every $\alpha\in\mathbb{R}^d$,
\begin{equation}\label{eq:minkowskiadd}
\textstyle h_{M(P)}(\alpha)=\sum_{k\in N_0}h_k(\alpha)\,,\qquad F_{M(P)}(\alpha)=\sum_{k\in N_0}F_{\Delta_k}(\alpha)\,.
\end{equation}
Moreover, $M(P)$ is full dimensional, because each $\Delta_k$ is.

\smallskip
\noindent\emph{Maximizing blocks.}
Recall the vertices $v^k_r$ ($r\in[d_0]$) of $\Delta_k$ from \cref{subsec:order}. For $k \in N_0$ and $\alpha \in \mathbb{R}^d$, define $s^k_r(\alpha) := \alpha^\top v^{k}_r = \sum_{t=r+1}^d \alpha_{\sigma^k_t}$\,, and let $R_k(\alpha) := \argmax_{r \in [d_0]}s^k_r(\alpha)$. We write the elements of $R_k(\alpha)$ in increasing order, $r^{k}_0 < r^{k}_1 < \cdots < r^{k}_{m_k}$\,, where $m_k = m_k(\alpha) := |R_k(\alpha)| - 1$. For $j = 1, \ldots, m_k$\,, the $j$-th \emph{maximizing block of $\sigma^k$ with respect
to $\alpha$} is $B_{k,j}(\alpha) := \{\sigma^k_{r^k_{j-1}+1}\,, \ldots, \sigma^k_{r^k_{j}} \}$, the set of entries of $\sigma^k$ in positions $r^k_{j-1}+1$ through $r^k_j$\,. Then $v^k_{r^k_{j-1}}-v^k_{r^k_j}=\mathbf{e}_{B_{k,j}(\alpha)}$\,. Entries before position $r^k_0+1$ or after position $r^k_{m_k}$ belong to no block; if $|R_k(\alpha)|=1$, there are no blocks, and if $\alpha=0$, every entry is a block by itself.

\begin{lemma}\label{lem:support}
For every $k\in N_0$ and $\alpha \in \mathbb{R}^d$:
\begin{enumerate}[label=(\alph*)]
    \item $h_k(\alpha)=\max_{r\in[d_0]}s_r^k(\alpha)$;
    \item $F_{\Delta_k}(\alpha)=\conv\{v_r^k:r\in R_k(\alpha)\}$, which is a simplex (a face of $\Delta_k$);
    \item $\dim F_{\Delta_k}(\alpha)=|R_k(\alpha)|-1=m_k(\alpha)$.
\end{enumerate}
\end{lemma}
\begin{proof}
By \cref{subsec:order}, the gap coordinates $y$ map $\Delta_k$ affinely onto $\Delta^{\mathrm{std}}_d$\,, sending $v^k_r$ to the $r$-th unit vector; so $\Delta_k$ is a $d$-simplex with vertices $v^k_0\,,\ldots,v^k_d$\,, and every $x\in\Delta_k$ is uniquely $x=\sum_{r}y_r v^k_r$ with $y\geq0$ and $\sum_r y_r=1$. With $h:=\max_r s^k_r(\alpha)$ we get $h-\alpha^\top x=\sum_r y_r\bigl(h-s^k_r(\alpha)\bigr)$, a sum of nonnegative terms that vanishes exactly when $y_r=0$ for all $r\notin R_k(\alpha)$. This gives (a) and (b), and (c) follows because the vertices $v^k_r$\,, $r\in R_k(\alpha)$, are affinely independent.
\end{proof}

Thus $h_k(\alpha)$ is the largest suffix sum of $\alpha$ along $\sigma^k$, including the empty suffix.

\begin{lemma}\label{lem:dimofFMP}
For every $k\in N_0$ and $\alpha\in\mathbb{R}^d$, the vectors $\mathbf{e}_{B_{k,j}(\alpha)}$\,, $1\leq j\leq m_k(\alpha)$, form a basis of $\dir F_{\Delta_k}(\alpha)$, and each of them is orthogonal to $\alpha$. Moreover, $\dim F_{M(P)}(\alpha) = \rank \{ \mathbf{e}_{B_{k,j}(\alpha)}: k \in N_0\,,\ 1 \leq j \leq m_k(\alpha) \}$. In particular, for $\alpha \neq 0$ this rank is at most $d-1$, and $\alpha$ is a facet normal of $M(P)$ if and only if the rank equals $d-1$, that is, if and only if the vectors $\mathbf{e}_{B_{k,j}(\alpha)}$ span $\alpha^\perp$.
\end{lemma}

\begin{proof}
By \cref{lem:support}, $F_{\Delta_k}(\alpha)=\conv\{v^k_r:r\in R_k(\alpha)\}$ is a simplex of dimension $m_k$ with vertices $v^k_{r^k_0}\,,\ldots,v^k_{r^k_{m_k}}$\,. Its direction space is spanned by the differences $v^k_{r^k_a}-v^k_{r^k_b}$ ($a<b$), and each such difference is the sum of the consecutive differences for $j=a+1,\ldots,b$. Hence $\dir F_{\Delta_k}(\alpha)$ is spanned by
\[
v^k_{r^{k}_{j-1}} - v^k_{r^{k}_{j}} = \textstyle\sum_{t=r^{k}_{j-1}+1}^{r^{k}_{j}} \mathbf{e}_{\sigma^k_t} = \mathbf{e}_{B_{k,j}(\alpha)}\,,\qquad 1\leq j\leq m_k\,.
\]
These $m_k$ vectors span an $m_k$-dimensional space, so they form a basis (equivalently, the blocks $B_{k,1}(\alpha),\ldots,B_{k,m_k}(\alpha)$ are nonempty and pairwise disjoint). Because $r^k_{j-1}\,,r^k_j\in R_k(\alpha)$, we have $\alpha^\top\mathbf{e}_{B_{k,j}(\alpha)}=s^k_{r^k_{j-1}}(\alpha)-s^k_{r^k_j}(\alpha)=0$.

Because $F_{M(P)}(\alpha)=\sum_{k\in N_0}F_{\Delta_k}(\alpha)$ by \eqref{eq:minkowskiadd} and direction spaces add, we get
\begin{align*}
\dir F_{M(P)}(\alpha) &= \textstyle\sum_{k \in N_0} \Span \{ \mathbf{e}_{B_{k,j}(\alpha)} : 1 \leq j \leq m_k \}\\
&= \Span \{ \mathbf{e}_{B_{k,j}(\alpha)} : k\in N_0\,,\ 1 \leq j \leq m_k \},
\end{align*}
which gives the dimension formula. For $\alpha\neq0$, this space lies in the $(d-1)$-dimensional space $\alpha^\perp$, so the rank is at most $d-1$. Finally, $M(P)$ is full dimensional, so $\alpha\neq0$ is a facet normal of $M(P)$ exactly when $\dim F_{M(P)}(\alpha)=d-1$, that is, when the block vectors span $\alpha^\perp$.
\end{proof}

\Cref{lem:dimofFMP} uses only that $\sigma^0,\ldots,\sigma^n$ are permutations of $V$, so it holds for every Minkowski sum of staircase simplices.

For nonempty $I\subset N_0$ put $Q_I:=\sum_{k\in I}\Delta_k$\,. Because the vertical facets of $\D$ correspond to the facet normals of $M(P)$ (\cref{sec:vertical}), we need that facet normals of $Q_I$ persist in $M(P)=Q_I+\sum_{k\in N_0\setminus I}\Delta_k$\,.

\begin{lemma}\label{lem:persistenceaddition}
Let $K \subset \mathbb{R}^d$ be a full-dimensional polytope, and let
$R \subset\mathbb{R}^d$ be a nonempty polytope. Every facet-normal ray of $K$ is a facet-normal ray of $K+R$. In particular, for every nonempty $I \subset N_0$\,, every facet-normal ray of $Q_I$ is a facet-normal ray of $M(P)$.
\end{lemma}

\begin{proof}
The polytope $K+R$ contains a translate of $K$, so it is full dimensional. If $\alpha$ is a facet normal of $K$ and $y\in F_R(\alpha)$, then $F_{K+R}(\alpha)=F_K(\alpha)+F_R(\alpha)$ (see \cite{Schneider2014}) contains the $(d-1)$-dimensional set $F_K(\alpha)+y$ and lies in a hyperplane, so it has dimension $d-1$. For the second statement, apply this to $K=Q_I$\,, which is full dimensional, and $R=\sum_{k\in N_0\setminus I}\Delta_k$ (if $I\neq N_0$).
\end{proof}

\cref{lem:persistenceaddition} reflects the familiar fact that the normal fan of $K+R$ is the common refinement of the normal fans of $K$ and $R$ (see \cite[Proposition~7.12]{ziegler}). The converse fails in general: $K+R$ may have facet normals that are facet normals of neither $K$ nor $R$, arising as new rays of the common refinement. By the equivalence of statements (2) and (3) of \cref{thm:verticalarelifting} below, $M(P)$ fails to be alcoved exactly when it has a facet normal that is a facet normal of no simplex of Stanley's triangulation.

\section{Vertical facets of \texorpdfstring{$\mathcal{D}$}{D} and full optimal big-M lifting}\label{sec:vertical}
We first show, for arbitrary full-dimensional polytopes, that the vertical facets of their Cayley embedding correspond to the facet normals of their Minkowski sum and are described by full optimal big-M liftings of supporting inequalities. We then specialize to Stanley's triangulation and prove \cref{thm:verticalarelifting}, which characterizes the posets for which lifting facet-describing inequalities of single simplices suffices.

\subsection{Full optimal big-M lifting and vertical facets}\label{subsec:lifted}
In this subsection, $\Po_0\,,\ldots,\Po_n\subset\mathbb{R}^d$ are arbitrary full-dimensional polytopes, $\D$ is their Cayley embedding as in \cref{sec:intro}, and we write points of $\mathbb{R}^{d+n}$ as $(x,z)$ with $x\in\mathbb{R}^d$ and $z\in\mathbb{R}^n$. Let $\bar\Po_k:=\Po_k\times\{\mathbf{e}_k\}$\,, so that $\D=\conv\bigl(\bigcup_{k\in N_0}\bar\Po_k\bigr)$, let $h_k:=h_{\Po_k}$\,, and let $M:=\sum_{k\in N_0}\Po_k$\,, which is full dimensional. By \cite[Theorem~2]{QU202531}, $\D$ is full dimensional, so the inequality describing a facet of $\D$ is unique up to positive scaling (\cref{subsec:support}).

\smallskip
\noindent\emph{Full optimal big-M lifting} (see \cite{QL_ISCO2024,QU202531}, and also \cite[Chapter II.2, Section 1]{NWbook}). Let
\begin{equation*}
\alpha^\top x \leq \beta \tag{\ensuremath{*}}\label{eq:star}
\end{equation*}
be valid for $\Po_k$ for some $k\in N_0$\,, and put $M_j := \min \{ \beta -\alpha^\top x : x \in \Po_j\}$ for $j\in N_0$\,. If $k=0$, lifting each $z_j$ with the largest coefficient that keeps the inequality valid on $\bar\Po_j$ gives
\begin{equation} \textstyle
\alpha^\top x + \sum_{i\in N} M_i z_i \leq \beta; \tag{\ensuremath{*_0}}\label{eq:star0}
\end{equation}
if $k\in N$, lifting $z_j$ in the same way for $j\in N\setminus\{k\}$, together with the term $M_0\bigl(1-\sum_{i\in N}z_i\bigr)$, which vanishes except on $\bar\Po_0$\,, gives
\begin{equation} \textstyle
\alpha^\top x + \sum_{i\in N \setminus\{k\}} (M_i-M_0) z_i -M_0 z_k\leq \beta  - M_0\,.
\tag{\ensuremath{*_k}}\label{eq:stark}
\end{equation}
Each coefficient affects validity on only one $\bar\Po_j$\,, so the order of lifting does not matter, and \eqref{eq:star0} and \eqref{eq:stark} are the \emph{full optimal big-M liftings} of \eqref{eq:star}. By \cite[Theorem~3]{QU202531}, the full optimal big-M lifting of a facet-describing inequality of some $\Po_k$ describes a facet of $\D$. The next result shows that lifting \emph{supporting} inequalities gives all vertical facets. For $\alpha \neq 0$, let
\begin{equation}\tag{\ensuremath{F_\alpha}}\label{eq:FP}
        \alpha^\top x + \textstyle\sum_{k \in N}(h_0(\alpha) - h_k(\alpha))z_k \leq h_0(\alpha).
\end{equation}

\begin{proposition}\label{prop:cayley}
Let $\alpha\neq0$.
\begin{enumerate}[label=(\alph*)]
    \item For every $k\in N_0$\,, \eqref{eq:FP} is the full optimal big-M lifting of the supporting inequality $\alpha^\top x\leq h_k(\alpha)$ of $\Po_k$\,; in particular, it is valid for $\D$.
    \item The face of $\D$ described by \eqref{eq:FP} is $\conv\bigl(\bigcup_{k\in N_0}(F_{\Po_k}(\alpha)\times\{\mathbf{e}_k\})\bigr)$, and its dimension is $n+\dim F_M(\alpha)$.
    \item \eqref{eq:FP} describes a facet of $\D$ if and only if $\alpha$ is a facet normal of $M$. Every vertical facet of $\D$ is described by \eqref{eq:FP} for some facet normal $\alpha$ of $M$, and $\{\lambda\alpha:\lambda>0\}\mapsto$ (the facet described by \eqref{eq:FP}) is a bijection between the facet-normal rays of $M$ and the vertical facets of $\D$.
\end{enumerate}
\end{proposition}

\begin{proof}
    (a) Take $\beta=h_k(\alpha)$ in \eqref{eq:star}. Then $M_j=h_k(\alpha)-h_j(\alpha)$ for $j\neq k$, and \eqref{eq:star0} (if $k=0$) and \eqref{eq:stark} (if $k\in N$, using $M_i-M_0=h_0(\alpha)-h_i(\alpha)$, $-M_0=h_0(\alpha)-h_k(\alpha)$, and $\beta-M_0=h_0(\alpha)$) both become \eqref{eq:FP}. At $(x,\mathbf{e}_k)$ with $x\in\Po_k$\,, \eqref{eq:FP} reads $\alpha^\top x\leq h_k(\alpha)$, which holds, with equality exactly for $x\in F_{\Po_k}(\alpha)$; so \eqref{eq:FP} is valid on every $\bar\Po_k$\,, hence on $\D$.

    (b) Let $\overline F$ be the face described by \eqref{eq:FP}. By (a), it contains every $(x,\mathbf{e}_k)$ with $x\in F_{\Po_k}(\alpha)$. Conversely, write $(\bar x,\bar z)\in\overline F$ as $\sum_k\lambda_k(x^k,\mathbf{e}_k)$ with $\lambda\geq0$, $\sum_k\lambda_k=1$, and $x^k\in\Po_k$\,. Then $\bar z_k=\lambda_k$ for $k\in N$, and the slack of \eqref{eq:FP} at $(\bar x,\bar z)$ is $\sum_k\lambda_k\bigl(h_k(\alpha)-\alpha^\top x^k\bigr)=0$, a sum of nonnegative terms; so $x^k\in F_{\Po_k}(\alpha)$ whenever $\lambda_k>0$. For the dimension, fix $\bar x^k\in F_{\Po_k}(\alpha)$ and put $w^k:=(\bar x^k-\bar x^0,\mathbf{e}_k)$ for $k\in N$. Every difference of generating points decomposes as $(x,\mathbf{e}_i)-(y,\mathbf{e}_j)=(x-\bar x^i,0)-(y-\bar x^j,0)+w^i-w^j$ (with $w^0:=0$), so $\dir\overline F=\bigl(U\times\{0\}\bigr)+\Span\{w^k:k\in N\}$ with $U:=\sum_k\dir F_{\Po_k}(\alpha)=\dir F_M(\alpha)$ by \eqref{eq:minkowskiadd}. The $z$-coordinates of the $w^k$ are linearly independent, so the sum is direct and $\dim\overline F=\dim F_M(\alpha)+n$. This is the dimension relation underlying the Cayley trick; see \cite{HuberRambauSantos2000}.

    (c) By (b), \eqref{eq:FP} describes a facet, $\dim\overline F=d+n-1$, if and only if $\dim F_M(\alpha)=d-1$, that is, $\alpha$ is a facet normal of $M$; the facet is vertical because $\alpha\neq0$. Now let $\alpha^\top x+\beta^\top z\leq\gamma$ describe a vertical facet $F$, so $\alpha\neq0$. If $F$ contained no point of $\bar\Po_k$\,, then, writing each point of $F$ as a convex combination of points of the $\bar\Po_j$ that are tight, $F$ would lie in the proper face $\D\cap\{z_k=0\}$ (if $k\in N$) or $\D\cap\{\sum_jz_j=1\}$ (if $k=0$), hence equal it, and $F$ would not be vertical. So for every $k$ some $x^k\in\Po_k$ satisfies $\alpha^\top x^k+\beta^\top\mathbf{e}_k=\gamma$; by validity, $x^k$ maximizes $\alpha^\top x$ over $\Po_k$\,, so $h_k(\alpha)+\beta^\top\mathbf{e}_k=\gamma$. With $\mathbf{e}_0=\mathbf{0}$ this gives $\gamma=h_0(\alpha)$ and $\beta_k=h_0(\alpha)-h_k(\alpha)$, so the inequality is \eqref{eq:FP}, and $\alpha$ is a facet normal of $M$ by the first statement. Finally, each $h_k$ is positively homogeneous, so the map is well defined on rays; it is surjective by what we just showed, and injective because two inequalities describing the same facet are positive multiples of each other, so their $x$-coefficients are too.
\end{proof}

\subsection{Specialization to Stanley's triangulation}\label{subsec:verticalnormals}
From now on, $\Po_k:=\Delta_k$\,, so that $\bar\Delta_k:=\Delta_k\times\{\mathbf{e}_k\}$ has extreme points $(v^k_r\,,\mathbf{e}_k)$, $r\in[d_0]$, $h_k$ is as in \cref{subsec:support}, and $M=M(P)$. \Cref{prop:cayley} then reads as follows.

\begin{corollary}\label{cor:bijection}
    For every facet normal $\alpha$ of $M(P)$, the inequality \eqref{eq:FP} describes a vertical facet of $\D$, and it is the full optimal big-M lifting of the supporting inequality $\alpha^\top x\leq h_k(\alpha)$ of \emph{every} $\Delta_k$\,; in particular, the choice of $\sigma^0$ affects only the description of the facet in indicator coordinates. The map sending a facet-normal ray $\{\lambda\alpha:\lambda>0\}$ of $M(P)$ to this facet is a bijection between the facet-normal rays of $M(P)$ and the vertical facets of $\D$. By \cref{lem:dimofFMP}, $\alpha\neq0$ is a facet normal of $M(P)$ if and only if the block vectors $\mathbf{e}_{B_{k,j}(\alpha)}$ have rank $d-1$.
\end{corollary}

\subsection{Vertical facets described by liftings}\label{subsec:whenlifting}
We now characterize the posets $P$ for which lifting facet-describing inequalities of single simplices suffices. We first isolate three lemmas.

\begin{lemma}\label{lem:liftedfacet}
Let $\alpha$ be a facet normal of $M(P)$, let $F$ be the vertical facet of $\D$ described by \eqref{eq:FP}, and let $k\in N_0$\,. Then $F$ is described by the full optimal big-M lifting of a facet-describing inequality of $\Delta_k$ if and only if $\alpha$ is a facet normal of $\Delta_k$\,.
\end{lemma}

\begin{proof}
If $\alpha$ is a facet normal of $\Delta_k$\,, then $\alpha^\top x\leq h_k(\alpha)$ is facet-describing for $\Delta_k$\,, and its full optimal big-M lifting is \eqref{eq:FP} by \cref{cor:bijection}. Conversely, suppose that $F$ is described by the full optimal big-M lifting of a facet-describing inequality $\tilde\alpha^\top x\leq\beta$ of $\Delta_k$\,, where $\tilde\alpha\in\mathbb{R}^d$ and $\beta\in\mathbb{R}$. Lifting does not change the $x$-coefficients, and the inequality describing $F$ is unique up to positive scaling, so $\tilde\alpha$ is a positive multiple of $\alpha$, which is therefore a facet normal of $\Delta_k$\,.
\end{proof}

\begin{lemma}\label{lem:alcovedtype}
Let $P$ be a poset on $V$, and let $i,j\in V$ be distinct. Then $-\mathbf{e}_i$ is a facet normal of $M(P)$ if and only if $i$ is minimal in $P$; $\mathbf{e}_i$ is a facet normal of $M(P)$ if and only if $i$ is maximal in $P$; and $\mathbf{e}_i-\mathbf{e}_j$ is a facet normal of $M(P)$ if and only if $i\parallel_P j$ or $i\prec_P j$ is a cover relation. In each case, the vector is a facet normal of $\Delta_\sigma$ for some $\sigma\in\mathcal L(P)$.
\end{lemma}

\begin{proof}
By \cref{lem:persistenceaddition}, a facet normal of some $\Delta_\sigma$\,, $\sigma\in\mathcal L(P)$, is a facet normal of $M(P)$. If $-\mathbf{e}_i$ is a facet normal of $M(P)$ and $j\prec_P i$ for some $j$, then every point of every $F_{\Delta_k}(-\mathbf{e}_i)$, and hence of $F_{M(P)}(-\mathbf{e}_i)=\sum_kF_{\Delta_k}(-\mathbf{e}_i)$, satisfies $x_i=x_j=0$, so this face has dimension at most $d-2$, a contradiction. Conversely, if $i$ is minimal, then $-\mathbf{e}_i$ is a facet normal of the simplex of a linear extension beginning with $i$. Symmetrically, $\mathbf{e}_i$ is a facet normal of $M(P)$ if and only if $i$ is maximal, and then it is a facet normal of the simplex of a linear extension ending with $i$.

Let $\alpha:=\mathbf{e}_i-\mathbf{e}_j$ be a facet normal of $M(P)$. If $j\prec_P i$, every maximizer of $\alpha$ on every $\Delta_k$ has $x_i=1$ and $x_j=0$, so $F_{M(P)}(\alpha)$ satisfies $x_i=n+1$ and $x_j=0$, a contradiction. If $i\prec_P t\prec_P j$ for some $t$, every maximizer has $x_i=x_t=x_j$\,, and these two independent equations again hold on $F_{M(P)}(\alpha)$. Hence $i\parallel_P j$ or $i\prec_P j$ is a cover relation.

Conversely, suppose that $i\parallel_P j$ or that $i\prec_P j$ is a cover relation; we construct a linear extension in which $i$ appears immediately before $j$. Let $I:=\{v\in V:v\prec_P i\text{ or }v\prec_P j\}\setminus\{i\}$. We first show that every predecessor of an element of $I$ also belongs to $I$. Suppose that $u\prec_P v$ for some $v\in I$. If $v\prec_P i$, then $u\prec_P i$, so $u\in I$. Otherwise, $v\prec_P j$, and hence $u\prec_P j$. Moreover, $u\ne i$, because $u=i$ would give $i\prec_P v\prec_P j$, contrary to the assumption. Thus $u\in I$.

    List the elements of $I$ in an order respecting $\prec_P$\,, then list $i,j$, and finally list the remaining elements in an order respecting $\prec_P$\,. Because every predecessor of $i$ lies in $I$ and every predecessor of $j$ in $I\cup\{i\}$, this is a linear extension of $P$, and $\mathbf{e}_i-\mathbf{e}_j$ is a facet normal of its simplex.
\end{proof}

\begin{lemma}\label{lem:compatfacts}
Let $\sigma,\tau$ be linear extensions of a poset on a finite ground set (such as $P$ or an induced subposet $P|_S$), with the element $0$ adjoined as in \eqref{eq:EKMdict}.
\begin{enumerate}
        \item[(A)]
        If deleting the same element $y$ from $\sigma$ and $\tau$
        gives the same ordering, then
        $\Delta_\sigma+\Delta_\tau$ is alcoved.

        \item[(B)]
        If there are three distinct elements $a,b,c$ such that
        $\sigma|_{\{a,b,c\}}=(a,b,c)$ and
        $\tau|_{\{a,b,c\}}=(c,b,a)$,
        then $\Delta_\sigma+\Delta_\tau$ is not alcoved.

        \item[(C)]
        If four distinct elements $a_1\,,a_2\,,a_3\,,a_4$ satisfy
        $\sigma|_{\{a_1,a_2,a_3,a_4\}}=(a_1\,,a_2\,,a_3\,,a_4)$ and
        $\tau|_{\{a_1,a_2,a_3,a_4\}}=(a_2\,,a_1\,,a_4\,,a_3)$,
        then $\Delta_\sigma+\Delta_\tau$ is not alcoved.
    \end{enumerate}
\end{lemma}

\begin{proof}
    By \eqref{eq:EKMdict}, $\Delta_\sigma+\Delta_\tau$ is alcoved if and only if the nondegenerate ordered set partitions $\mathbf S_\sigma=(\sigma_d\,,\ldots,\sigma_1\,,0)$ and $\mathbf S_\tau$ (written as sequences of singletons) are compatible. By \cite[Theorem~5.5]{early2025when}, two nondegenerate ordered set partitions are compatible if and only if they are neither $4$- nor $6$-interlaced in the sense of \cite[Definition~5.4]{early2025when}. By \cite[Corollary~4.5]{early2025when}, the normal fan of $\Delta_{\mathbf S}$ is invariant under cyclic shifts of $\mathbf S$; because whether a Minkowski sum is alcoved depends only on the normal fans of the summands, replacing $\mathbf S_\sigma$ or $\mathbf S_\tau$ by a cyclic shift does not affect compatibility.

    \emph{Proof of (A).} A nondegenerate ordered set partition, written as a sequence of letters, determines a cyclic order; by the preceding paragraph, compatibility depends only on the two cyclic orders, and accordingly the patterns of \cite[Definition~5.4]{early2025when} are matched up to cyclic rotation. For a cyclic order, write $x\to x'$ if $x'$ immediately follows $x$. Because $\sigma$ and $\tau$ agree after deleting $y$, the cyclic orders of $\mathbf S_\sigma$ and $\mathbf S_\tau$ agree after deleting $y$; hence, for every subset $A$ of the ground set, the cyclic orders of $\mathbf S_\sigma|_A$ and $\mathbf S_\tau|_A$ coincide if $y\notin A$, and coincide after deleting $y$ if $y\in A$. It therefore suffices to show that, in each pattern pair of \cite[Definition~5.4]{early2025when}, the two cyclic orders remain different after deleting any single letter (and without deleting one):
    \begin{itemize}
        \item $(w_1\,,w_2\,,w_3\,,w_4)$ and $(w_3\,,w_2\,,w_1\,,w_4)$: the successor of $w_1$ is $w_2$\,, resp.\ $w_4$; after deleting $w_2$\,, $w_3$\,, or $w_4$ it is $w_3$ versus $w_4$\,, $w_2$ versus $w_4$\,, $w_2$ versus $w_3$; after deleting $w_1$\,, the successor of $w_4$ is $w_2$ versus $w_3$\,.
        \item $(w_1\,,\ldots,w_6)$ and $(w_3\,,w_4\,,w_1\,,w_2\,,w_5\,,w_6)$: the successor of $w_2$ is $w_3$\,, resp.\ $w_5$\,, and after deleting any letter other than $w_2$ the two successors of $w_2$ are still different; after deleting $w_2$\,, the successor of $w_1$ is $w_3$ versus $w_5$\,.
        \item $(w_1\,,\ldots,w_6)$ and $(w_1\,,w_4\,,w_5\,,w_2\,,w_3\,,w_6)$: the successor of $w_1$ is $w_2$\,, resp.\ $w_4$\,, and after deleting any letter other than $w_1$ the two successors of $w_1$ are still different; after deleting $w_1$\,, the successor of $w_6$ is $w_2$ versus $w_4$\,.
    \end{itemize}
    Therefore $\mathbf S_\sigma$ and $\mathbf S_\tau$ are neither $4$- nor $6$-interlaced, and they are compatible.

    \emph{Proof of (B).} The restrictions of $\mathbf S_\sigma$ and $\mathbf S_\tau$ to $\{a,b,c,0\}$ are $(c,b,a,0)$ and $(a,b,c,0)$; writing $(w_1\,,w_2\,,w_3\,,w_4):=(c,b,a,0)$, these are $(w_1\,,w_2\,,w_3\,,w_4)$ and $(w_3\,,w_2\,,w_1\,,w_4)$, so $\mathbf S_\sigma$ and $\mathbf S_\tau$ are $4$-interlaced and not compatible.

    \emph{Proof of (C).} Write $A:=\{a_1\,,a_2\,,a_3\,,a_4\}$. We have $\mathbf S_\sigma|_{A}=(a_4\,,a_3\,,a_2\,,a_1)$ and $\mathbf S_\tau|_{A\cup\{0\}}=(a_3\,,a_4\,,a_1\,,a_2\,,0)$. Let $\mathbf T$ be the cyclic shift of $\mathbf S_\tau$ that begins with $a_2$; then $\mathbf T|_{A}=(a_2\,,a_3\,,a_4\,,a_1)$. With $(w_1\,,w_2\,,w_3\,,w_4):=(a_4\,,a_3\,,a_2\,,a_1)$ these are $(w_1\,,w_2\,,w_3\,,w_4)$ and $(w_3\,,w_2\,,w_1\,,w_4)$, so $\mathbf S_\sigma$ and $\mathbf T$ are $4$-interlaced. Hence $\mathbf S_\sigma$ and $\mathbf T$, and therefore $\mathbf S_\sigma$ and $\mathbf S_\tau$\,, are not compatible.
\end{proof}

\begin{theorem}[Vertical facets and full optimal big-M liftings]\label{thm:verticalarelifting}
Let $P$ be a poset on $V=[d]$ with $\mathcal L(P)=\{\sigma^0,\ldots,\sigma^n\}$. The following statements are equivalent:
\begin{enumerate}
    \item Every vertical facet of $\D$ is described by the full optimal big-M lifting of a facet-describing inequality of some simplex $\Delta_{k}$\,, $k \in N_0$\,.
    \item Every facet normal of $M(P)$ is a facet normal of some simplex $\Delta_{k}$\,, $k \in N_0$\,.
    \item $M(P)$ is an alcoved polytope.
    \item The simplices $\Delta_{k}$\,, $k \in N_0$\,, are pairwise compatible, that is, $\Delta_{k_1}+\Delta_{k_2}$ is alcoved for all $k_1\,,k_2 \in N_0$ with $k_1 \neq k_2$\,.
    \item The graph $G(P)$ is acyclic.
    \item There exists $y \in V$ such that $V \setminus \{y\}$ is a chain.
\end{enumerate}
In particular, if $P$ does not satisfy (6), then $\D$ has a vertical facet that is not described by the full optimal big-M lifting of any facet-describing inequality of any $\Delta_k$\,.
\end{theorem}

\begin{proof}
    We prove $(1)\Leftrightarrow(2)\Leftrightarrow(3)\Leftrightarrow(4)\Leftrightarrow(6)\Leftrightarrow(5)$.

    \noindent\emph{$(1)\Leftrightarrow(2)$.} By \cref{cor:bijection}, the vertical facets of $\D$ are exactly the facets described by \eqref{eq:FP} for the facet normals $\alpha$ of $M(P)$, so this is \cref{lem:liftedfacet}.

    \noindent\emph{$(2)\Rightarrow(3)$.} By \cref{subsec:order}, every facet normal of $\Delta_k$ is a positive multiple of $-\mathbf{e}_{\sigma^k_1}$\,, $\mathbf{e}_{\sigma^k_d}$\,, or $\mathbf{e}_{\sigma^k_r}-\mathbf{e}_{\sigma^k_{r+1}}$ ($1\leq r<d$). So under (2), every facet normal of $M(P)$ has the form required in \cref{subsec:alcoved}.

    \noindent\emph{$(3)\Rightarrow(2)$.} By (3), every facet normal of $M(P)$ is a positive multiple of some $-\mathbf{e}_i$\,, $\mathbf{e}_i$\,, or $\mathbf{e}_i-\mathbf{e}_j$\,, and by \cref{lem:alcovedtype} it is a facet normal of some $\Delta_k$\,.

    \noindent\emph{$(3)\Leftrightarrow(4)$.} If $n=0$, both statements hold. Otherwise, \cref{lem:persistenceaddition} shows that every facet-normal ray of any pairwise sum $\Delta_{k_1}+\Delta_{k_2}$ is a facet-normal ray of $M(P)$. Hence $(3)$ implies $(4)$. The converse is the pairwise compatibility theorem \cite[Theorem~A]{early2025when}, transferred to our coordinates by $\psi$ (\cref{subsec:alcoved}).

    \noindent\emph{$(6)\Rightarrow(4)$.} Suppose that (6) holds, and write $C=V\setminus\{y\}=\{c_1\,,\ldots,c_m\}$ with $c_1\prec_P\cdots\prec_P c_m$\,, where $m=d-1$. Let $\sigma,\tau\in\mathcal L(P)$. Both linear extensions respect the order on $C$, so deleting $y$ from each gives $\sigma|_C=\tau|_C=(c_1\,,\ldots,c_m)$. By \cref{lem:compatfacts}(A), $\Delta_\sigma+\Delta_\tau$ is alcoved. Hence the simplices of Stanley's triangulation are pairwise compatible, proving (4).

    \noindent\emph{$(4)\Rightarrow(6)$.} Suppose that (4) holds. By \cref{lem:restriction}(3), for every nonempty $S\subset V$, $\Delta_\rho+\Delta_\eta$ is alcoved for all $\rho,\eta\in\mathcal L(P|_S)$. Suppose, for a contradiction, that (6) fails. By \cref{lem:posetchain}(c), one of the following two cases occurs.
    If $a,b,c$ are pairwise incomparable, put $S=\{a,b,c\}$; then $\rho=(a,b,c)$ and $\eta=(c,b,a)$ are linear extensions of $P|_S$\,, and by \cref{lem:compatfacts}(B) $\Delta_\rho+\Delta_\eta$ is not alcoved.
    Otherwise there are distinct $b_1\,,b_2\,,t_1\,,t_2$ with $b_1\parallel_P b_2$\,, $t_1\parallel_P t_2$\,, and no $t_i$ below a $b_j$\,; put $S:=\{b_1\,,b_2\,,t_1\,,t_2\}$. Then $\rho=(b_1\,,b_2\,,t_1\,,t_2)$ and $\eta=(b_2\,,b_1\,,t_2\,,t_1)$ are linear extensions of $P|_S$\,, and by \cref{lem:compatfacts}(C) $\Delta_\rho+\Delta_\eta$ is not alcoved.
    Both cases contradict (4), so (6) holds.

    \noindent\emph{$(6)\Rightarrow(5)$.} Suppose that (6) holds, and write $C=V\setminus\{y\}=\{c_1\,,\ldots,c_m\}$ with $c_1\prec_P\cdots\prec_P c_m$\,, where $m=d-1$. By transitivity, the elements of $C$ below $y$ form an initial segment, and those above $y$ form a final segment. Thus, there exist $0\leq p\leq q\leq m$ such that
    \[
        \begin{aligned}
        c_i&\prec_P y &&(i\leq p),\\
        c_i&\parallel_P y &&(p<i\leq q),\\
        y&\prec_P c_i &&(i>q).
        \end{aligned}
    \]
    The linear extensions of $P$ are precisely $\tau^t=(c_1\,,\ldots,c_t\,,y,c_{t+1}\,,\ldots,c_m)$, $p\leq t\leq q$. By \cref{lem:adjacentswap}, every edge of $G(P)$ interchanges two incomparable elements in adjacent positions; because every two elements of $C$ are comparable, such an interchange involves $y$ and changes its position by one. Hence $E(G(P))=\bigl\{\{\tau^t,\tau^{t+1}\}:p\leq t<q\bigr\}$, and $G(P)$ is a path (a single vertex if $p=q$). In particular, $G(P)$ is acyclic, proving (5).

    \noindent\emph{$(5)\Rightarrow(6)$.} Suppose that (5) holds, and suppose, for a contradiction, that (6) fails. By \cref{lem:posetchain}(c), one of the following two cases occurs.

    First suppose that $a,b,c$ are pairwise incomparable, and let $I:=\{v\in V:v\prec_P a\text{ or }v\prec_P b\text{ or }v\prec_P c\}$. If $v\in I$ and $u\prec_P v$, then transitivity gives $u\in I$. List the elements of $I$ in an order respecting $\prec_P$; then all predecessors of $a,b,c$ have been listed, and because $a,b,c$ are pairwise incomparable they may be listed next in any order; finally, list the remaining elements in an order respecting $\prec_P$\,, using the same orderings before and after $a,b,c$ in all six cases. This gives six linear extensions of $P$, which we consider in the order $(a,b,c)$, $(b,a,c)$, $(b,c,a)$, $(c,b,a)$, $(c,a,b)$, $(a,c,b)$, displaying only the positions occupied by $a,b,c$. Each consecutive pair, including the last and the first, differs by interchanging two incomparable elements in adjacent positions, so by \cref{lem:adjacentswap} these six linear extensions form a six-cycle in $G(P)$, contradicting (5).

    Otherwise there are distinct $b_1\,,b_2\,,t_1\,,t_2$ with $b_1\parallel_P b_2$\,, $t_1\parallel_P t_2$\,, and no $t_i$ below a $b_j$\,. Let $I:=\{v\in V: v\prec_P b_1\text{ or }v\prec_P b_2\}$ and $J:=I\cup\{b_1\,,b_2\}\cup\{v\in V: v\prec_P t_1\text{ or }v\prec_P t_2\}$. Both $I\cup\{b_1\,,b_2\}$ and $J$ are closed under taking predecessors, $b_1\,,b_2\notin I$, and $t_1\,,t_2\notin J$ (because $t_1\parallel_P t_2$ and no $t_i$ lies below a $b_j$). List $I$ in an order respecting $\prec_P$\,, then $b_1\,,b_2$\,, then $J\setminus(I\cup\{b_1\,,b_2\})$ in an order respecting $\prec_P$\,, then $t_1\,,t_2$\,, and finally the remaining elements in an order respecting $\prec_P$\,. This is a linear extension of $P$ in which $b_1\,,b_2$ are adjacent and $t_1\,,t_2$ are adjacent. Interchanging $b_1\,,b_2$\,, interchanging $t_1\,,t_2$\,, or both gives four linear extensions which, by \cref{lem:adjacentswap}, form a four-cycle in $G(P)$, contradicting (5). Hence (6) holds.
\end{proof}

\begin{figure}[htbp]
\centering
\begin{tikzpicture}[scale=0.75,
    x=1cm, y=1cm,
    vertex/.style={
        circle, fill=black,
        inner sep=0pt,
        minimum size=3.5pt,
        outer sep=0pt
    },
    every label/.append style={font=\small},
    line width=0.5pt
]

\begin{scope}
    \node[vertex,label=above:{$a_1$}] at (0,0.65) {};
    \node[vertex,label=above:{$a_2$}] at (0.625,0.65) {};
    \node[vertex,label=above:{$a_3$}] at (1.25,0.65) {};

    \node at (0.625,-1.2) {$A_3$};
\end{scope}

\begin{scope}[xshift=2.9cm]
    \node[vertex,label=below:{$b_1$}] (bL1) at (0,0) {};
    \node[vertex,label=below:{$b_2$}] (bL2) at (1.25,0) {};
    \node[vertex,label=above:{$t_1$}] (bU1) at (0,1.3) {};
    \node[vertex,label=above:{$t_2$}] (bU2) at (1.25,1.3) {};

    \draw (bL1)--(bU1)
          (bL2)--(bU2);

    \node at (0.625,-1.2) {$2+2$};
\end{scope}

\begin{scope}[xshift=5.8cm]
    \node[vertex,label=below:{$b_1$}] (cL1) at (0,0) {};
    \node[vertex,label=below:{$b_2$}] (cL2) at (1.25,0) {};
    \node[vertex,label=above:{$t_1$}] (cU1) at (0,1.3) {};
    \node[vertex,label=above:{$t_2$}] (cU2) at (1.25,1.3) {};

    \draw (cL1)--(cU1)
          (cL2)--(cU1)
          (cL2)--(cU2);

    \node at (0.625,-1.2) {$\mathsf{N}$};
\end{scope}

\begin{scope}[xshift=8.7cm]
    \node[vertex,label=below:{$b_1$}] (dL1) at (0,0) {};
    \node[vertex,label=below:{$b_2$}] (dL2) at (1.25,0) {};
    \node[vertex,label=above:{$t_1$}] (dU1) at (0,1.3) {};
    \node[vertex,label=above:{$t_2$}] (dU2) at (1.25,1.3) {};

    \draw (dL1)--(dU1)
          (dL1)--(dU2)
          (dL2)--(dU1)
          (dL2)--(dU2);

    \node at (0.625,-1.2) {$K_{2,2}$};
\end{scope}

\end{tikzpicture}
\caption{The four minimal forbidden induced subposets of \cref{cor:forbidden}.}
\label{fig:heighttwoobstructions}
\end{figure}
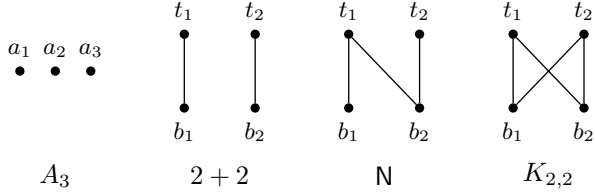

Let $A_3$ be the three-element antichain, and let $2+2$, $\mathsf{N}$, and $K_{2,2}$ be the posets on $\{b_1\,,b_2\,,t_1\,,t_2\}$ whose relations are, respectively, $\{b_1\prec t_1\,,\ b_2\prec t_2\}$, $\{b_1\prec t_1\,,\ b_2\prec t_1\,,\ b_2\prec t_2\}$, and $\{b_i\prec t_j: i,j\in\{1,2\}\}$; see \cref{fig:heighttwoobstructions}.

\begin{corollary}\label{cor:forbidden}
    Let $P$ be a poset on $V=[d]$. The equivalent conditions (1)--(6) of \cref{thm:verticalarelifting} hold if and only if $P$ contains no induced subposet isomorphic to $A_3$\,, $2+2$, $\mathsf{N}$, or $K_{2,2}$\,. These four posets are the minimal such obstructions: none of them satisfies (6), and every proper induced subposet of each of them does. Equivalently, all edges of the incomparability graph of $P$ share a common vertex.
\end{corollary}

\begin{proof}
    Condition (6) is inherited by induced subposets: if $V\setminus\{y\}$ is a chain, then so is $S\setminus\{y\}$ for every $S\subset V$. None of the four posets satisfies (6): in $A_3$ any two elements are incomparable, and each of $2+2$, $\mathsf{N}$, $K_{2,2}$ contains the two disjoint incomparable pairs $\{b_1\,,b_2\}$ and $\{t_1\,,t_2\}$\,, so deleting one element always leaves an incomparable pair. Hence (6) implies that $P$ contains none of them.

    Conversely, suppose that $P$ contains none of the four posets, and suppose, for a contradiction, that (6) fails. Because $P$ has no three-element antichain, \cref{lem:posetchain}(c) gives distinct $b_1\,,b_2\,,t_1\,,t_2$ with $b_1\parallel_P b_2$\,, $t_1\parallel_P t_2$\,, and no $t_i$ below a $b_j$\,; put $S:=\{b_1\,,b_2\,,t_1\,,t_2\}$. Every relation in $P|_S$ has the form $b_i\prec_P t_j$\,. No element of $S$ is incomparable to the other three, because together with the two elements of the other pair it would form an $A_3$\,. Thus the relations $b_i\prec_P t_j$ form a bipartite graph on $\{b_1\,,b_2\}\sqcup\{t_1\,,t_2\}$ without isolated vertices, with $2$, $3$, or $4$ edges, and $P|_S$ is isomorphic to $2+2$, $\mathsf{N}$, or $K_{2,2}$\,, respectively, a contradiction. Hence (6) holds.

    For minimality, every proper induced subposet of $A_3$ has at most two elements and satisfies (6) trivially. A proper induced subposet of $2+2$, $\mathsf{N}$, or $K_{2,2}$ with three elements contains one of the comparable pairs $b_1\prec t_1$\,, $b_2\prec t_2$\,, so deleting the third element leaves a chain; smaller ones satisfy (6) trivially. Finally, all incomparable pairs contain a common element $y$ exactly when $V\setminus\{y\}$ is a chain, which is (6).
\end{proof}

\Cref{lem:alcovedtype,lem:liftedfacet} identify the vertical facets of $\D$ that are described by liftings of facet-describing inequalities of single simplices, and they count the vertical facets of $\D$ whenever \cref{thm:verticalarelifting} applies. For a poset $P$, let $\min(P)$ and $\max(P)$ denote the numbers of minimal and maximal elements, $\operatorname{inc}(P)$ the number of unordered incomparable pairs, and $\operatorname{cov}(P)$ the number of cover relations.

\begin{corollary}\label{cor:alcovedlist}
Let $P$ be a poset on $V=[d]$. The vertical facets of $\D$ that are described by full optimal big-M liftings of facet-describing inequalities of simplices of Stanley's triangulation are exactly those described by \eqref{eq:FP} for the vectors $\alpha$ listed in \cref{lem:alcovedtype}, and their number is $\min(P)+\max(P)+2\operatorname{inc}(P)+\operatorname{cov}(P)$. In particular, under the conditions of \cref{thm:verticalarelifting}, this is the number of vertical facets of $\D$.
\end{corollary}

\begin{proof}
The facet normals of a simplex $\Delta_\sigma$ are, up to positive scaling, $-\mathbf{e}_{\sigma_1}$\,, $\mathbf{e}_{\sigma_d}$\,, and $\mathbf{e}_{\sigma_r}-\mathbf{e}_{\sigma_{r+1}}$ (\cref{subsec:order}); by \cref{lem:alcovedtype}, as $\sigma$ ranges over $\mathcal L(P)$, these are exactly the listed vectors. By \cref{cor:bijection,lem:liftedfacet}, a vertical facet of $\D$ is described by the full optimal big-M lifting of a facet-describing inequality of some $\Delta_k$ if and only if its normal is a facet normal of some $\Delta_k$\,. The listed vectors span pairwise distinct rays, and there are $\min(P)+\max(P)+2\operatorname{inc}(P)+\operatorname{cov}(P)$ of them, so the count follows from \cref{cor:bijection}. Under the conditions of \cref{thm:verticalarelifting}, statement (2) says that $M(P)$ has no other facet normals.
\end{proof}

For a chain, \cref{cor:alcovedlist} gives $1+1+0+(d-1)=d+1$, the number of facets of the simplex $\D=\Delta_0$\,. For the antichain on two elements it gives $2+2+2+0=6$: $M(P)$ is the sum of the two triangles $\Delta_{(1,2)}$ and $\Delta_{(2,1)}$\,, a hexagon.

\begin{example}\label[example]{ex:antichain}
    The smallest poset violating (6) is the antichain $P$ on $V=[3]$. All six permutations are linear extensions; we order them as $\sigma^0=(1,2,3)$, $\sigma^1=(1,3,2)$, $\sigma^2=(2,1,3)$, $\sigma^3=(2,3,1)$, $\sigma^4=(3,1,2)$, $\sigma^5=(3,2,1)$. By \cref{lem:adjacentswap}, $G(P)$ is the six-cycle $(1,2,3)$, $(2,1,3)$, $(2,3,1)$, $(3,2,1)$, $(3,1,2)$, $(1,3,2)$. Let $\alpha:=\mathbf{e}_1+\mathbf{e}_2-\mathbf{e}_3$\,. Its suffix sums along the six permutations give $(h_0(\alpha),\ldots,h_5(\alpha))=(1,1,1,1,2,2)$, and the only maximizing blocks are $B=\{1,3\}$ for $\sigma^1=(1,3,2)$ and $B=\{2,3\}$ for $\sigma^3=(2,3,1)$. The block vectors $(1,0,1)$ and $(0,1,1)$ have rank $2=d-1$, so $\alpha$ is a facet normal of $M(P)$ by \cref{lem:dimofFMP}, and by \cref{cor:bijection} the inequality $x_1+x_2-x_3-z_4-z_5\leq 1$ describes a vertical facet of $\D$. Because $\alpha$ is not a positive multiple of any $\pm\mathbf{e}_i$ or $\mathbf{e}_i-\mathbf{e}_j$\,, it is not a facet normal of any $\Delta_k$\,, so this facet is not described by the full optimal big-M lifting of a facet-describing inequality of any $\Delta_k$\,. Only the two simplices $\Delta_1$ and $\Delta_3$ contribute blocks, so $\alpha$ is already a facet normal of $\Delta_1+\Delta_3$; this sum of two staircase simplices is \SF{\SFexample}.
\end{example}

\begin{remark}[Consequences of condition (6)]\label[remark]{rem:height2}
    Condition (6) of \cref{thm:verticalarelifting} also bounds the number of simplices. If $V\setminus\{y\}$ is a chain $c_1\prec_P\cdots\prec_P c_{d-1}$\,, then the chain elements below $y$ form an initial segment and those above $y$ a final segment, so the elements incomparable to $y$ are consecutive, say $c_{i+1}\,,\ldots,c_{i+k}$\,, and the linear extensions of $P$ are obtained by inserting $y$ into this segment. Hence $|\mathcal L(P)|=k+1\leq d$, with equality exactly when $y$ is incomparable to every other element. Thus the full optimal big-M liftings of the facet-describing inequalities of the simplices, together with the bounds on $z$, describe $\D$ only when Stanley's triangulation has at most $d$ simplices, whereas in general $|\mathcal L(P)|$ can be as large as $d!$.

    Condition (6) requires a chain with $d-1$ elements, so it fails whenever $\height(P)\leq d-2$. In particular, if $\height(P)\leq2$, the conditions of \cref{thm:verticalarelifting} hold if and only if $d\leq2$, or $d=3$ and $P$ is not an antichain; equivalently, $d\leq3$ and $\width(P)\leq2$. Thus for every poset of height at most two with $d\geq4$, the graph $G(P)$ contains a cycle, $M(P)$ is not alcoved, and $\D$ has a vertical facet that is not described by the full optimal big-M lifting of a facet-describing inequality of any $\Delta_k$\,. In terms of \cref{cor:forbidden}, one checks directly that a poset of height at most two contains $A_3$ if and only if $\width(P)\geq3$, and it contains $2+2$, $\mathsf N$, or $K_{2,2}$ if and only if it has two relations $b_1\prec_P t_1$ and $b_2\prec_P t_2$ on four distinct elements, because for such elements the absence of three-element chains leaves $b_1\prec_P t_2$ and $b_2\prec_P t_1$ as the only possible further relations.
\end{remark}

\section{Posets close to a chain}\label{sec:nearchain}
We now use the order relations of $P$ to bound the entries of the primitive facet normals of $M(P)$. \Cref{subsec:crossproduct} expresses facet normals as generalized cross products. In \cref{subsec:generalbound} we show that they are bounded by a polynomial in $d$ of degree $\ell(P)-1$, and that they lie in $\{-1,0,1\}$ when $\ell(P)\leq1$. In \cref{subsec:chaintwo} we treat $\ell(P)=2$: we determine all facet normals of $M(P)$ for a chain with two isolated elements, and deduce the sharp bound $d-2$; by \cref{cor:bijection}, this also describes all vertical facets of $\D$ for that poset. In \cref{subsec:lowerbound} we show that this degree is attained: for every fixed $\ell\geq2$ there are posets with $\ell(P)=\ell$ whose facet normals have entries of order $d^{\ell-1}$, and, when $\ell$ grows with $d$, entries exponential in $d$. In \cref{subsec:upper3} we prove the upper bound $\lfloor(d-1)^2/2\rfloor$ for $\ell(P)=3$, make the leading constants for $\ell(P)=4$ explicit, and formulate the remaining open problems (\cref{rem:open}).

\subsection{Facet normals as generalized cross products}\label{subsec:crossproduct}
By \cref{lem:dimofFMP}, a facet normal $\alpha$ of $M(P)$ is determined, up to scaling, by $d-1$ maximizing blocks of linear extensions of $P$ whose indicator vectors are linearly independent. The next lemma makes this explicit; it is the basic tool for bounding the entries of $\alpha$. In the companion paper \cite{QL_staircase} it is used for Minkowski sums of arbitrary staircase simplices.

\begin{lemma}[Facet normals as generalized cross products]\label{lem:crossproduct}
    Let $B_1\,,\ldots,B_{d-1}\subset[d]$ be nonempty subsets whose indicator vectors are linearly independent. Let $A$ be the $(d-1)\times d$ matrix with rows $\mathbf{e}_{B_j}^\top$\,, let $\widehat A^i$ be obtained from $A$ by deleting column $i$, and let $\gamma_i:=(-1)^{i-1}\det\widehat A^i$ for $i\in[d]$ (the determinant of the empty matrix being one). We call $\gamma(B_1\,,\ldots,B_{d-1}):=\gamma$ the \emph{generalized cross product} of $\mathbf{e}_{B_1}\,,\ldots,\mathbf{e}_{B_{d-1}}$\,. Then $\gamma\in\mathbb{Z}^d\setminus\{0\}$ and $\ker A=\Span\{\gamma\}$. Consequently, if $\alpha\in\mathbb{R}^d\setminus\{0\}$ satisfies $\sum_{i\in B_j}\alpha_i=0$ for all $j\in[d-1]$, then $\alpha=t\gamma$ for some $t\neq0$; if moreover $\alpha$ is a primitive integer vector, then $\alpha=\pm\gamma/g$ with $g:=\gcd(\gamma_1\,,\ldots,\gamma_d)$, and in particular $|\alpha_i|\leq|\det\widehat A^i|$ for every $i\in[d]$. \end{lemma}

\begin{proof}
    The vector $\gamma$ is integral because $A$ is a zero-one matrix. For each $j\in[d-1]$, let $\widetilde A_j$ be the $d\times d$ matrix whose first row is row $j$ of $A$ and whose remaining rows are the rows of $A$. It has two identical rows, so expanding along the first row gives $0=\det\widetilde A_j=\sum_{i=1}^{d}(-1)^{1+i}A_{ji}\det\widehat A^i=(A\gamma)_j\,$. Hence $A\gamma=0$. Because $\rank A=d-1$, some maximal minor $\det\widehat A^i$ is nonzero, so $\gamma\neq0$; and $\ker A$ is one-dimensional, so $\ker A=\Span\{\gamma\}$.

    If $\alpha\neq0$ satisfies the hypothesis, then $A\alpha=0$, so $\alpha\in\ker A$ and $\alpha=t\gamma$ with $t\neq0$. Suppose moreover that $\alpha$ is a primitive integer vector. Then $\gamma=t^{-1}\alpha$, and B\'ezout's identity gives integers $\lambda_i$ with $\sum_i\lambda_i\alpha_i=1$, so $t^{-1}=\sum_i\lambda_i\gamma_i$ is a nonzero integer. Thus $\gamma$ is a nonzero integer multiple of the primitive vector $\alpha$, which forces $t^{-1}=\pm g$, that is, $\alpha=\pm\gamma/g$. In particular $|\alpha_i|\leq|\gamma_i|=|\det\widehat A^i|$.
\end{proof}

Each maximizing block of a linear extension $\sigma$ of $P$ is a run of consecutive entries of $\sigma$. Consequently, every maximizing block $B$ is \emph{order-convex} in $P$: if $a\prec_P b\prec_P c$ and $a,c\in B$, then $b\in B$.

For a nonnegative integer $j$, let $\delta_j$ be the greatest absolute value of the determinant of a $j\times j$ zero-one matrix, with $\delta_0:=1$; this is the zero-one version of Hadamard's maximal determinant problem (see \cite{OEISA003432}). Bordering a matrix by a corner entry $1$ shows that $\delta_j$ is nondecreasing in $j$, and expanding along a row gives $\delta_j\leq j\,\delta_{j-1}$\,, so $\delta_j\leq j!$. The first values are $\delta_1=\delta_2=1$, $\delta_3=2$, $\delta_4=3$, and $\delta_5=5$ (see \cite{OEISA003432}).

\begin{lemma}\label{lem:delta4}
For every $j\geq1$, $\delta_j\leq(j+1)^{(j+1)/2}/2^j$ \textup{(}see \cite[Lemma~24]{Ziegler01}\textup{)}. In particular, $\delta_4=3$.
\end{lemma}

\begin{proof}
The bound is \cite[Lemma~24]{Ziegler01}. For $j=4$ it gives $\delta_4\leq5^{5/2}/16<4$, so $\delta_4\leq3$; conversely, with $J_4$ the all-ones matrix, $J_4-I_4$ has eigenvalues $3,-1,-1,-1$, so $|\det(J_4-I_4)|=3$.
\end{proof}

For $j=5$, \cref{lem:delta4} gives only $\delta_5\leq6$; the exact value $\delta_5=5$ is taken from \cite{OEISA003432}.

\subsection{A bound in terms of the distance to a chain}\label{subsec:generalbound}
Let $\ell(P):=d-\height(P)$ be the least number of elements whose deletion from $P$ leaves a chain, and for integers $0\leq\ell<d$ let
\[
H_\ell(d):=\max\bigl\{\lVert\alpha\rVert_\infty:\ \alpha \text{ a primitive facet normal of } M(P)\bigr\},
\]
the maximum being over all posets $P$ on $[d]$ with $\ell(P)=\ell$. Recall the maximal zero-one determinants $\delta_j$ from \cref{subsec:crossproduct}. Throughout this subsection, $C=\{c_1\,,\ldots,c_m\}$ with $c_1\prec_P\cdots\prec_P c_m$ is a chain of $P$, $E:=V\setminus C$, and $\ell:=|E|=d-m$; from \cref{thm:posetheight} on, $C$ is a longest chain, so that $\ell=\ell(P)$.

\begin{lemma}\label{lem:posetbounds}
Every primitive facet normal $\alpha$ of $M(P)$ satisfies
\begin{align*}
 |\alpha_i|&\leq\tbinom{d-1}{\ell}\delta_\ell
 &&\text{for } i\in C,\\
 |\alpha_i|&\leq\tbinom{d-1}{\ell-1}\delta_{\ell-1}
 &&\text{for } i\in E \text{ (if } \ell\geq1).
\end{align*}
\end{lemma}

\begin{proof}
Let $\alpha$ be a primitive facet normal of $M(P)$. By \cref{lem:dimofFMP}, we may choose $d-1$ maximizing blocks $B_1\,,\ldots,B_{d-1}$ of the linear extensions of $P$ such that $\mathbf{e}_{B_1}\,,\ldots,\mathbf{e}_{B_{d-1}}$ are linearly independent. Let $A$ be the $(d-1)\times d$ matrix whose $j$-th row is $\mathbf{e}_{B_j}^{\top}$. Because $\sum_{i\in B_j}\alpha_i=0$ for $j\in[d-1]$, we have $A\alpha=0$ and $\rank A=d-1$.

\emph{The chain columns are totally unimodular.} Let $A_C$ be the submatrix of $A$ formed by the columns indexed by $C$, in the order $c_1\,,\ldots,c_m$\,. Every linear extension of $P$ lists the elements of $C$ in this order, and each $B_j$ consists of consecutive entries of a linear extension (so it is order-convex, \cref{subsec:crossproduct}). Hence, if $c_a\,,c_b\in B_j$ with $a<b$, then $c_t\in B_j$ for $a<t<b$: the ones in each row of $A_C$ are consecutive. Matrices with this consecutive-ones property (interval matrices) are totally unimodular (see \cite{NWbook}).

\emph{Reduction to minors.} For $i\in[d]$, let $\widehat A^{\,i}$ be obtained from $A$ by deleting column $i$. By \cref{lem:crossproduct}, $|\alpha_i|\leq|\det\widehat A^{\,i}|$ for every $i\in[d]$. For row and column index sets $I$ and $J$, write $A_{I,J}$ for the corresponding submatrix of $A$; the determinant of the $0\times0$ matrix is one.

\emph{The bounds.} Let $i\in C$. The matrix $\widehat A^{\,i}$ contains all $\ell$ columns indexed by $E$, and the generalized Laplace expansion along these columns gives
\[
\bigl|\det\widehat A^{\,i}\bigr|
\leq
\textstyle\sum_{I\subset[d-1],\,|I|=\ell}
\bigl|\det A_{I,E}\bigr|\,
\bigl|\det A_{[d-1]\setminus I,\;C\setminus\{i\}}\bigr|.
\]
Each first factor is a zero-one determinant of order $\ell$, so it is at most $\delta_\ell$; each second factor uses only columns indexed by $C$, so it is at most $1$; and there are $\binom{d-1}{\ell}$ choices of $I$. This gives the bound for $i\in C$. If $i\in E$, then $\widehat A^{\,i}$ contains the $\ell-1$ columns indexed by $E\setminus\{i\}$, and the same argument applies with $\ell-1$ in place of $\ell$.

\end{proof}

We now bound the chain coordinates of a primitive facet normal by a multiple of its coordinates on $E$, which lowers the degree in $d$ of the resulting bound from $\ell$ to $\ell-1$. We use prefix sums. For a linear extension $\sigma^k$ and $\alpha\in\mathbb{R}^d$, let $\Phi^k_r:=\sum_{t=1}^{r}\alpha_{\sigma^k_t}$ for $0\leq r\leq d$. Because $\alpha^\top v^k_r=\sum_{i\in[d]}\alpha_i-\Phi^k_r$\,, a position $r$ belongs to $R_k(\alpha)$ if and only if $\Phi^k_r$ attains its minimum over $0\leq r\leq d$; we call such a position \emph{minimal}. Thus a maximizing block is a run of $\sigma^k$ between two consecutive minimal positions.

Put $\Sigma_i:=\sum_{j\leq i}\alpha_{c_j}$ for $0\leq i\leq m$ (so $\Sigma_0=0$) and $\Sigma_{\min}:=\min_{0\leq i\leq m}\Sigma_i$\,. Every linear extension of $P$ lists the elements of $C$ in the order $c_1\,,\ldots,c_m$\,, and a maximizing block $B$ is a run of a linear extension, so $B\cap C=\{c_{r+1}\,,\ldots,c_s\}$, where $r$ and $s$ are the numbers of elements of $C$ before $B$ and through the end of $B$; we call $r$ and $s$ the \emph{endpoint indices} of $B$. The two endpoints of $B$ are minimal positions at which exactly $r$, respectively $s$, elements of $C$ have appeared. Finally, for $v\in\mathbb{R}^E$ put $\lVert v\rVert_\infty:=\max_{e\in E}|v_e|$\,, $W_+(v):=\sum_{e\in E}\max(v_e\,,0)$, $W_-(v):=\sum_{e\in E}\max(-v_e\,,0)$, and $W(v):=\max\bigl(W_+(v),W_-(v)\bigr)=\max_{S\subset E}\bigl|\sum_{e\in S}v_e\bigr|\leq\ell\,\lVert v\rVert_\infty$\,, and write $\alpha_E$ for the restriction of $\alpha$ to $E$.

\begin{lemma}\label{lem:endpoints}
Let $\alpha$ be a primitive facet normal of $M(P)$ with $\lVert\alpha\rVert_\infty\geq2$. Then every $i\in\{0,1,\ldots,m\}$ is an endpoint index of some maximizing block.
\end{lemma}

\begin{proof}
Suppose that some $i$ is an endpoint index of no maximizing block. A maximizing block with endpoint indices $r,s$ contains $c_1$ only if $r=0$; it contains $c_m$ only if $s=m$; and it contains exactly one of $c_i\,,c_{i+1}$ (for $0<i<m$) only if $r=i$ or $s=i$. Hence one of $\mathbf{e}_{c_1}$ (if $i=0$), $\mathbf{e}_{c_i}-\mathbf{e}_{c_{i+1}}$ (if $0<i<m$), or $\mathbf{e}_{c_m}$ (if $i=m$) is orthogonal to the indicator vectors of all maximizing blocks. By \cref{lem:dimofFMP}, these indicator vectors span $\alpha^\perp$, so that vector is a nonzero multiple of $\alpha$; because $\alpha$ is a primitive integer vector and the entries of the vector lie in $\{-1,0,1\}$, we get $\lVert\alpha\rVert_\infty=1$, a contradiction.
\end{proof}

\begin{lemma}[Window lemma]\label{lem:window}
Let $\alpha\in\mathbb{R}^d$, and let $\lambda$ satisfy $\Sigma_\lambda=\Sigma_{\min}$\,. If $i$ is an endpoint index of a maximizing block of some linear extension of $P$, then
\[
\Sigma_i\leq\Sigma_{\min}+W_+(\alpha_E)\ \ \text{if } i\leq\lambda,
\qquad
\Sigma_i\leq\Sigma_{\min}+W_-(\alpha_E)\ \ \text{if } i\geq\lambda.
\]
If $\alpha$ is a primitive facet normal of $M(P)$ with $\lVert\alpha\rVert_\infty\geq2$, then consequently $\Sigma_{\min}\leq\Sigma_i\leq\Sigma_{\min}+W(\alpha_E)$ for $0\leq i\leq m$, and $|\alpha_c|\leq W(\alpha_E)\leq\ell\,\lVert\alpha_E\rVert_\infty$ for every $c\in C$.
\end{lemma}

\begin{proof}
Let $i$ be an endpoint index of a maximizing block of $\sigma^k$\,; if $i=\lambda$ there is nothing to prove. Then $\sigma^k$ has a minimal position $p$ at which exactly $i$ elements of $C$ have appeared. Because the elements of $C$ appear in $\sigma^k$ in their chain order, the position $q$ at which exactly $\lambda$ elements of $C$ have appeared, taken as the position of $c_\lambda$ (the start of $\sigma^k$ if $\lambda=0$), satisfies $q>p$ if $i<\lambda$ and $q<p$ if $i>\lambda$. Let $S\subset E$ be the set of elements of $E$ in positions $\min(p,q)+1,\ldots,\max(p,q)$ of $\sigma^k$, and put $\alpha(S):=\sum_{e\in S}\alpha_e$\,.

If $i<\lambda$, the entries of $\sigma^k$ in positions $p+1,\ldots,q$ are $c_{i+1}\,,\ldots,c_\lambda$ and the elements of $S$, so $\Phi^k_q-\Phi^k_p=(\Sigma_{\min}-\Sigma_i)+\alpha(S)$; minimality of $p$ gives $\Phi^k_q-\Phi^k_p\geq0$, so $\Sigma_i\leq\Sigma_{\min}+\alpha(S)\leq\Sigma_{\min}+W_+(\alpha_E)$. If $i>\lambda$, the entries in positions $q+1,\ldots,p$ are $c_{\lambda+1}\,,\ldots,c_i$ and the elements of $S$, so $\Phi^k_p-\Phi^k_q=(\Sigma_i-\Sigma_{\min})+\alpha(S)$; minimality of $p$ gives $\Phi^k_p-\Phi^k_q\leq0$, so $\Sigma_i\leq\Sigma_{\min}-\alpha(S)\leq\Sigma_{\min}+W_-(\alpha_E)$.

For the second statement, every $i$ is an endpoint index by \cref{lem:endpoints}, and $W_\pm(\alpha_E)\leq W(\alpha_E)$. Finally, $\alpha_{c_i}=\Sigma_i-\Sigma_{i-1}$ with both terms in $[\Sigma_{\min}\,,\Sigma_{\min}+W(\alpha_E)]$, so $|\alpha_{c_i}|\leq W(\alpha_E)$; and $W(\alpha_E)\leq\ell\,\lVert\alpha_E\rVert_\infty$ because $W_+(\alpha_E)$ and $W_-(\alpha_E)$ are sums of at most $\ell$ terms, each at most $\lVert\alpha_E\rVert_\infty$\,.
\end{proof}

\begin{theorem}\label{thm:posetheight}
    Let $\ell(P)\geq1$. Every primitive facet normal $\alpha$ of $M(P)$ satisfies
    \[
    \lVert\alpha\rVert_\infty\,\leq\,\delta_{\ell(P)}\tbinom{d-1}{\ell(P)-1}\,\leq\,\ell(P)\,(d-1)^{\ell(P)-1}.
    \]
    In particular, for posets with bounded $\ell(P)=d-\height(P)$, the entries of the primitive facet normals of $M(P)$ grow only polynomially in $d$, with degree $\ell(P)-1$: for each fixed $\ell\geq1$, $H_\ell(d)\leq(\delta_\ell/(\ell-1)!)\,d^{\ell-1}+O(d^{\ell-2})$.
\end{theorem}

\begin{proof}
    Here $C$ is a longest chain, so $\ell=\ell(P)\geq1$. If $\lVert\alpha\rVert_\infty=1$, the first bound holds because both factors are at least one. Otherwise, by \cref{lem:dimofFMP} choose $d-1$ maximizing blocks $B_1\,,\ldots,B_{d-1}$ of linear extensions of $P$ whose indicator vectors are linearly independent, let $A$ be the $(d-1)\times d$ matrix with rows $\mathbf{e}_{B_j}^\top$\,, and let $\gamma$ be its cofactor vector; by \cref{lem:crossproduct}, $\alpha=\pm\gamma/g$ with $g:=\gcd(\gamma_1\,,\ldots,\gamma_d)\geq1$.

    Fix a nonempty $S\subset E$, and let $A_S$ be the $d\times d$ matrix obtained from $A$ by appending the row $\mathbf{e}_S^\top$\,. In the expansion of $\det A_S$ along its last row, the cofactor of the entry in column $i$ is $(-1)^{d+i}\det\widehat A^i=(-1)^{d+1}\gamma_i$\,, so $\det A_S=(-1)^{d+1}\sum_{e\in S}\gamma_e$\,. Now expand $\det A_S$ by the generalized Laplace expansion along the $d-\ell$ columns indexed by $C$. The last row of $A_S$ vanishes on these columns, so only the $\binom{d-1}{d-\ell}=\binom{d-1}{\ell-1}$ choices of $d-\ell$ rows among the first $d-1$ rows contribute. Each such minor in the chain columns has absolute value at most one, because its rows have consecutive ones in the chain order and interval matrices are totally unimodular (proof of \cref{lem:posetbounds}). Each complementary minor is the determinant of an $\ell\times\ell$ zero-one matrix, so it has absolute value at most $\delta_\ell$\,. Hence $\bigl|\sum_{e\in S}\gamma_e\bigr|\leq\delta_\ell\binom{d-1}{\ell-1}$ for every $S\subset E$, and therefore $W(\alpha_E)=\max_{S\subset E}|\sum_{e\in S}\alpha_e|=\max_{S\subset E}|\sum_{e\in S}\gamma_e|/g\leq\delta_\ell\binom{d-1}{\ell-1}$. By \cref{lem:window}, $|\alpha_c|\leq W(\alpha_E)$ for every $c\in C$, and $|\alpha_e|\leq W(\alpha_E)$ for every $e\in E$ by taking $S=\{e\}$. This proves the first inequality.

    For the second, $\delta_\ell\leq\ell!$ (\cref{subsec:crossproduct}) and $\binom{d-1}{\ell-1}(\ell-1)!\leq(d-1)^{\ell-1}$. The asymptotic statement follows from $\binom{d-1}{\ell-1}=d^{\ell-1}/(\ell-1)!+O(d^{\ell-2})$.
\end{proof}

\begin{remark}\label[remark]{rem:posetheight}
Combining the bound of \cref{lem:posetbounds} for $i\in E$ with \cref{lem:window} gives $\lVert\alpha\rVert_\infty\leq\ell\,\delta_{\ell-1}\binom{d-1}{\ell-1}$; \cref{thm:posetheight} replaces $\ell\,\delta_{\ell-1}$ by $\delta_\ell$\,, which is never larger, because $\delta_\ell\leq\ell\,\delta_{\ell-1}$\,. For $\ell=2,3,4,5$ this factor drops from $2,3,8,15$ to $1,2,3,5$. At the other extreme, $C=\varnothing$ in \cref{lem:posetbounds} gives the bound $\delta_{d-1}$\,, which holds for every Minkowski sum of staircase simplices (see \cite{QL_staircase}). For small $\ell$, the bound of \cref{thm:posetheight} is $1$ for $\ell=1$, $d-1$ for $\ell=2$, $(d-1)(d-2)$ for $\ell=3$, and $(d-1)(d-2)(d-3)/2$ for $\ell=4$, and sharper results are available for $\ell\leq2$: if $\ell\leq1$, then $V\setminus\{y\}$ is a chain for some $y\in V$, so $M(P)$ is alcoved by \cref{thm:verticalarelifting} and every facet normal is a positive multiple of some $\pm\mathbf{e}_i$ or $\mathbf{e}_i-\mathbf{e}_j$\,; hence $\lVert\alpha\rVert_\infty=1$, matching the bound. If $\ell=2$, then $\lVert\alpha\rVert_\infty\leq d-2$, and this is sharp, by \cref{thm:sharpchaintwoelements} below. Thus \cref{thm:posetheight} gives new information only for $\ell\geq3$. In \cref{subsec:lowerbound} we show that its degree $\ell-1$ is attained: for every fixed $\ell\geq2$ there are posets with entries of order $d^{\ell-1}$ (\cref{cor:growth}), and in \cref{subsec:upper3} we improve the constant for $\ell=3$ and make it explicit for $\ell=4$; for $\ell=3$ the largest entry lies between $\lfloor(d-2)^2/3\rfloor-1$ and $\lfloor(d-1)^2/2\rfloor$ (\cref{cor:grid,thm:upper3}); the exact values remain open (\cref{rem:open}). When $\ell(P)$ grows with $d$, the entries can be exponential in $d$ (\cref{cor:exp}); by contrast, sums of arbitrary staircase simplices can have facet normals with entries exponential in $d$ already for three summands, and of order $d^{d/2}2^{-O(d)}$ when there are many summands (see \cite{QL_staircase}).
\end{remark}

We close this subsection with a lemma, valid for every poset, that removes coordinates on which $\alpha$ vanishes; it is used in \cref{subsec:chaintwo}.

\begin{lemma}\label{lem:deletezeros}
    Let $P$ be a poset on $V$, let $\alpha\in\mathbb{R}^V\setminus\{0\}$, and let $Z\subseteq V$ satisfy $\alpha_z=0$ for every $z\in Z$. Set $\overline{P}:=P|_{V\setminus Z}$, and let $\overline{\alpha}$ be the restriction of $\alpha$ to $V\setminus Z$. If $\alpha$ is a facet normal of $M(P)$, then $\overline{\alpha}$ is a facet normal of $M(\overline{P})$. Conversely, suppose that, for every $z\in Z$, the singleton $\{z\}$ is a maximizing block of some linear extension of $P$ with respect to $\alpha$. If $\overline{\alpha}$ is a facet normal of $M(\overline{P})$, then $\alpha$ is a facet normal of $M(P)$.
\end{lemma}

\begin{proof}
    Let $\pi$ be the coordinate projection from $\mathbb{R}^V$ onto $\mathbb{R}^{V\setminus Z}$, and put $U:=\dir F_{M(P)}(\alpha)$ and $\overline U:=\dir F_{M(\overline P)}(\overline\alpha)$. For $\sigma\in\mathcal L(P)$, \cref{lem:restriction} gives $\pi(\Delta_\sigma)=\Delta_{\sigma|_{V\setminus Z}}$\,, and every linear extension of $\overline P$ arises as such a restriction. Because $\alpha$ vanishes on $Z$, we have $\alpha^\top x=\overline\alpha^\top\pi(x)$; hence $\pi$ maps $F_{\Delta_\sigma}(\alpha)$ onto $F_{\Delta_{\sigma|_{V\setminus Z}}}(\overline\alpha)$, and $\pi(\dir F_{\Delta_\sigma}(\alpha))=\dir F_{\Delta_{\sigma|_{V\setminus Z}}}(\overline\alpha)$. Summing over $\sigma\in\mathcal L(P)$ and using \eqref{eq:minkowskiadd} gives $\pi(U)=\overline U$. By hypothesis, $\mathbf{e}_z\in U$ for every $z\in Z$, so $\ker\pi=\Span\{\mathbf{e}_z:z\in Z\}\subset U$, and therefore $\dim U=|Z|+\dim\overline U$. Thus $\dim U=|V|-1$ if and only if $\dim\overline U=|V\setminus Z|-1$; because $\overline\alpha\neq0$, \cref{lem:dimofFMP} completes the proof of the first statement. For the second, if $\alpha$ is a facet normal, then $U$ has dimension $|V|-1$ and lies in $\alpha^\perp$, so $U=\alpha^\perp$, which contains $\mathbf{e}_z$ for every $z\in Z$; the same argument then gives $\dim\overline U=|V\setminus Z|-1$.
\end{proof}

For Minkowski sums of arbitrary staircase simplices, deleting zero coordinates does not change the least number of summands needed to realize a facet normal \SF{\SFzeropadding}; no hypothesis on singleton blocks is needed there, because the summands need not be linear extensions of a fixed poset.

\subsection{A chain and two elements: a complete classification}\label{subsec:chaintwo}
For $m\geq1$ and $\ell\geq0$, let $P^{(m,\ell)}$ be a chain $c_1\prec\cdots\prec c_m$ together with $\ell$ isolated elements, which are incomparable with each other and with every chain element. Every poset on the same ground set in which $c_1\prec\cdots\prec c_m$ is a chain is obtained by adding relations to $P^{(m,\ell)}$, so its linear extensions form a subset of $\mathcal L(P^{(m,\ell)})$, and by \cref{lem:persistenceaddition} its facet normals are facet normals of $M(P^{(m,\ell)})$. In particular, after relabeling a longest chain, the maximum $H_\ell(d)$ of \cref{subsec:generalbound} is attained by $P^{(d-\ell,\ell)}$\,.

In this subsection $\ell=2$. We write $P^{(m)}:=P^{(m,2)}$, $C_m:=\{c_1\,,\ldots,c_m\}$, and $a,b$ for the two isolated elements; thus $d=m+2$, and we index coordinates by $c_1\,,\ldots,c_m\,,a,b$. Because a linear extension of $P^{(m)}$ is determined by the positions of $a$ and $b$, there are $(m+1)(m+2)$ of them. We use the prefix sums $\Phi^k_r$ and minimal positions of \cref{subsec:generalbound}.

In \cref{lem:twoisolatessigns,lem:twoisolatesreduce} and \cref{thm:posettwoisolates}, for $\alpha\in\mathbb{R}^d$ with $\alpha_a\,,\alpha_b>0$ we write $u:=\alpha_a$\,, $v:=\alpha_b$\,, $L:=u+v$, $\omega:=-\sum_{i=1}^m\alpha_{c_i}$\,, and $S:=\{c\in C_m:\alpha_c<0\}$, $Z:=C_m\setminus S$.

\begin{lemma}\label{lem:twoisolatessigns}
    Let $\alpha$ be a facet normal of $M(P^{(m)})$ with $\alpha_a\,,\alpha_b>0$. Then $\alpha_c\leq0$ for every $c\in C_m$\,, $S\neq\varnothing$, and $\omega\leq L$.
\end{lemma}

\begin{proof}
    By \cref{lem:dimofFMP}, the vectors orthogonal to all block indicators are the multiples of $\alpha$; as $\alpha_a\neq0$, none of them is a nonzero vector supported on $C_m$\,.

    Suppose that $\alpha_{c_i}>0$ for some $i<m$, and fix a linear extension $\sigma^k$ in which $c_i$ and $c_{i+1}$ occupy positions $\rho_1<\rho_2$\,. Only $a$ and $b$ can lie strictly between them, and their coefficients are positive, so $\Phi^k_r>\Phi^k_{\rho_1-1}\geq\min\Phi^k$ for $\rho_1\leq r<\rho_2$\,. So no position in $[\rho_1\,,\rho_2)$ is minimal, no block contains exactly one of $c_i\,,c_{i+1}$\,, and $\mathbf{e}_{c_i}-\mathbf{e}_{c_{i+1}}$ is orthogonal to all block indicators, a contradiction. If $\alpha_{c_m}>0$, no position at or after that of $c_m$ is minimal, so $\mathbf{e}_{c_m}$ is orthogonal to all block indicators. Hence $\alpha_c\leq0$ on $C_m$\,.

    If $S=\varnothing$, a block, having $\alpha$-sum zero, cannot contain $a$, so $\mathbf{e}_a$ is orthogonal to all block indicators. Finally, suppose that $\omega>L$, and let $s$ be the first element of $S$ in the chain order. In every linear extension, the prefix sums before $s$ are nonnegative, whereas the total sum is $L-\omega<0$; hence no position before $s$ is minimal, and $\mathbf{e}_s$ is orthogonal to all block indicators.
\end{proof}

\begin{lemma}\label{lem:twoisolatesreduce}
    Let $\alpha\in\mathbb{R}^d$ satisfy $\alpha_a\,,\alpha_b>0$, $\alpha_c\leq0$ for every $c\in C_m$\,, $S\neq\varnothing$, and $\omega\leq L$. Then $\alpha$ is a facet normal of $M(P^{(m)})$ if and only if its restriction $\overline\alpha$ to $S\cup\{a,b\}$ is a facet normal of $M(\overline P)$, where $\overline P:=P^{(m)}|_{S\cup\{a,b\}}$ is a chain with two isolated elements.
\end{lemma}

\begin{proof}
    By \cref{lem:deletezeros}, it suffices to show that $\{z\}$ is a maximizing block for every $z\in Z$. In a linear extension in which $a$ and $b$ come immediately after $z$, the prefix sums decrease to some $-t$ just before $z$, equal $-t$ just after $z$, increase by $L$ at $a$ and $b$, and then decrease to $L-\omega\geq0$. So $\{z\}$ is a maximizing block.
\end{proof}

The next lemma treats the reduced case, with all chain coefficients negative.

\begin{lemma}\label{lem:twoisolatescycle}
    Let $h\geq1$, and let $\overline P$ be a chain $s_1\prec\cdots\prec s_h$ together with two isolated elements $a,b$. Let $\beta$ have $\beta_a=u$, $\beta_b=v$ and $\beta_{s_i}=-w_i$\,, where $u,v,w_1\,,\ldots,w_h$ are positive integers, put $t_0:=0$ and $t_i:=w_1+\cdots+w_i$\,, and suppose that $\omega:=t_h\leq L:=u+v$. Let $g:=\gcd(u,v)$, $p:=u/g$ and $q:=v/g$. Then $\beta$ is a facet normal of $M(\overline P)$ if and only if $w_1=\cdots=w_h=g$ and $h\in\{p+q-1,\,p+q\}$.
\end{lemma}

\begin{proof}
    \emph{The maximizing blocks.} A block has $\beta$-sum zero and all chain coefficients are negative, so it contains $a$ or $b$, and its chain elements form an interval $\{s_{i+1}\,,\ldots,s_j\}$\,. If a block contains $a$ but not $b$, then $t_j-t_i=u$; if it contains $b$ but not $a$, then $t_j-t_i=v$; and a block containing both has chain weight $L$, which, because $\omega\leq L$ and all $w_i>0$, happens only when $\omega=L$ and the block is the whole ground set. Conversely, if $t_j-t_i=u$, then along the linear extension $(s_1\,,\ldots,s_i\,,\ a,\ s_{i+1}\,,\ldots,s_j\,,\ b,\ s_{j+1}\,,\ldots,s_h)$ the prefix sums have minimum $-t_i$ on the first chain portion, $u-t_j=-t_i$ on the second, and $L-\omega\geq0$ on the third, and they exceed $-t_i$ strictly inside $\{a,s_{i+1}\,,\ldots,s_j\}$; so this set is a maximizing block. Interchanging $a$ and $b$ treats $t_j-t_i=v$. If $\omega=L$, the linear extension $(a,b,s_1\,,\ldots,s_h)$ has minimal prefix sum only at its two ends, so the whole ground set is a maximizing block. Hence the block indicators are exactly $\mathbf{e}_a+\sum_{k=i+1}^{j}\mathbf{e}_{s_k}$ for $t_j-t_i=u$, $\mathbf{e}_b+\sum_{k=i+1}^{j}\mathbf{e}_{s_k}$ for $t_j-t_i=v$, and, when $\omega=L$, $\mathbf{e}_a+\mathbf{e}_b+\sum_{k=1}^{h}\mathbf{e}_{s_k}$\,.

    \emph{The orthogonal complement.} By \cref{lem:dimofFMP}, $\beta$ is a facet normal if and only if the space $K$ of vectors orthogonal to all these indicators is one-dimensional. For $y\in\mathbb{R}^{S\cup\{a,b\}}$, put $x_0:=0$ and $x_i:=-\sum_{k=1}^{i}y_{s_k}$ for $i\in[h]$; this is a linear bijection between $y$ and $(y_a\,,y_b\,,x_1\,,\ldots,x_h)$. Orthogonality to these indicators reads $x_j-x_i=y_a$ whenever $t_j-t_i=u$, $x_j-x_i=y_b$ whenever $t_j-t_i=v$, and, when $\omega=L$, $x_h=y_a+y_b$\,.

    On the integers $\{0,\ldots,L-1\}$, consider the directed graph with arcs $t\to t+u\pmod L$. It consists of $g$ disjoint directed cycles, one for each residue class modulo $g$, each with $p+q$ vertices. Let $\mathcal G$ be the subgraph induced by $\mathcal Q:=\{t_0\,,\ldots,t_h\}$ if $\omega<L$, and by $\mathcal Q:=\{t_0\,,\ldots,t_{h-1}\}$ if $\omega=L$; in the latter case we identify $t_h=L$ with $0$ and use $x_h=y_a+y_b$ to eliminate $x_h$\,. Writing $x(t_i):=x_i$\,, the equations become
    \begin{equation}\label{eq:xtu}
        x(t+u)-x(t)=y_a\quad\text{if }t+u<L,\qquad
        x(t-v)-x(t)=-y_b\quad\text{if }t+u\geq L,
    \end{equation}
    one for each arc of $\mathcal G$, because $t+u\equiv t-v\pmod L$; at the identified endpoint, the equations involving $L$ reduce to the same form after substituting $x_h=y_a+y_b$\,.

    Every connected component of $\mathcal G$ is a directed path (possibly a single vertex) or a complete directed cycle. On a path, \eqref{eq:xtu} determines all potentials from one initial potential and $y_a\,,y_b$\,. On a complete cycle, let $n_u$ and $n_v$ be the numbers of arcs of the first and second form in \eqref{eq:xtu}. Summing the numerical increments around the cycle gives $n_uu-n_vv=0$ and $n_u+n_v=p+q$, so $n_u=q$ and $n_v=p$, because $\gcd(p,q)=1$; summing the potential increments shows that the equations on the cycle are consistent exactly when $qy_a-py_b=0$. Let $\nu$ be the number of components of $\mathcal G$, and let $\varepsilon:=1$ if $\mathcal G$ contains a complete cycle and $\varepsilon:=0$ otherwise. All complete cycles impose the same equation, and because $x(0)=0$, every component other than the one containing $0$ has one free initial potential. Therefore $\dim K=(\nu-1)+(2-\varepsilon)=\nu+1-\varepsilon$.
    Thus $\dim K=1$ if and only if $\mathcal G$ is connected and contains a complete cycle, that is, if and only if $\mathcal Q$ is the cycle $\{0,g,2g,\ldots,L-g\}$ through $0$. Because $0=t_0<t_1<\cdots$, this means $t_i=ig$ for every $i$, with $t_h=L-g$ (so $h=p+q-1$) if $\omega<L$, and $t_h=L$ (so $h=p+q$) if $\omega=L$. Equivalently, $w_1=\cdots=w_h=g$ and $h\in\{p+q-1,p+q\}$.
\end{proof}

\begin{theorem}\label{thm:posettwoisolates}
The primitive facet normals $\alpha$ of $M(P^{(m)})$ satisfying $\alpha_a>0$ and $\alpha_b>0$ are exactly the vectors $\alpha=-\mathbf{e}_S+p\,\mathbf{e}_a+q\,\mathbf{e}_b$\,, where $\varnothing\neq S\subset C_m$\,, $p,q\in\mathbb{Z}_{>0}$\,, $\gcd(p,q)=1$, and $|S|\in\{p+q-1,p+q\}$. The primitive facet normals satisfying $\alpha_a<0$ and $\alpha_b<0$ are exactly the negatives of these vectors.
\end{theorem}

\begin{proof}
    \emph{Necessity.} Let $\alpha$ be a primitive facet normal such that $\alpha_a>0$ and $\alpha_b>0$. By \cref{lem:twoisolatessigns}, $\alpha_c\leq0$ on $C_m$\,, $S\neq\varnothing$ and $\omega\leq L$, so by \cref{lem:twoisolatesreduce} the restriction $\overline\alpha$ to $S\cup\{a,b\}$ is a facet normal of $M(\overline P)$, where $\overline P$ is the chain $S$ together with $a,b$. Because $\alpha$ is integral, \cref{lem:twoisolatescycle} applies with $h:=|S|$ and gives $\alpha_c=-g$ for $c\in S$ and $|S|\in\{p+q-1,p+q\}$, where $g=\gcd(u,v)$, $p=u/g$, $q=v/g$. Hence $\alpha=g(-\mathbf{e}_S+p\,\mathbf{e}_a+q\,\mathbf{e}_b)$, and because $\alpha$ is primitive and $S\neq\varnothing$, we get $g=1$.

    \emph{Sufficiency.} Let $\alpha=-\mathbf{e}_S+p\,\mathbf{e}_a+q\,\mathbf{e}_b$ with $S$, $p$, $q$ as in the statement. Its chain coefficients lie in $\{0,-1\}$, $S\neq\varnothing$, and $\omega=|S|\leq p+q=L$, so \cref{lem:twoisolatesreduce} applies. The restriction $\overline\alpha$ satisfies the hypotheses of \cref{lem:twoisolatescycle} with $u=p$, $v=q$, $g=1$, $w_i=1$ and $h=|S|\in\{p+q-1,p+q\}$, so it is a facet normal of $M(\overline P)$; hence $\alpha$ is a facet normal of $M(P^{(m)})$. It is primitive, because its entries on $S$ equal $-1$.

    \emph{Negative coordinates at $a$ and $b$.} Let $R$ reverse the chain coordinates ($R\mathbf{e}_{c_i}=\mathbf{e}_{c_{m+1-i}}$) and fix $\mathbf{e}_a$\,, $\mathbf{e}_b$\,. Reversing a linear extension and relabeling $c_i\leftrightarrow c_{m+1-i}$ gives a linear extension, so $x\mapsto\mathbf{e}-Rx$ permutes the simplices of $\mathcal{T}(P^{(m)})$ and sends facet normals $\alpha$ of $M(P^{(m)})$ to $-R\alpha$. As $R$ preserves $|S|$, this shows that the primitive facet normals with $\alpha_a<0$ and $\alpha_b<0$ are exactly the negatives of the vectors above.
\end{proof}

In \cref{thm:posettwoisolates} only the size of $S$ matters; $S$ need not be an interval of the chain. Because $|S|\leq m$, we have $p+q\leq m+1$, so $\lVert\alpha\rVert_\infty=\max\{p,q\}\leq m=d-2$.

\begin{corollary}\label{cor:posetsharp}
    Let $m\geq1$. The vector $\alpha^*:=-\mathbf{e}_{C_m}+\mathbf{e}_a+m\,\mathbf{e}_b$ is a primitive facet normal of $M(P^{(m)})$ with $\lVert\alpha^*\rVert_\infty=m=d-2$. Moreover, among the primitive facet normals $\alpha$ of $M(P^{(m)})$ with $\alpha_a\alpha_b>0$, equality $\lVert\alpha\rVert_\infty=m$ holds exactly for $\alpha=\pm\alpha^*$ and $\alpha=\pm(-\mathbf{e}_{C_m}+m\,\mathbf{e}_a+\mathbf{e}_b)$.
\end{corollary}

\begin{proof}
    With $S=C_m$\,, $p=1$ and $q=m$, we have $\gcd(p,q)=1$ and $|S|=m=p+q-1$, so $\alpha^*$ is a primitive facet normal by \cref{thm:posettwoisolates}, and clearly $\lVert\alpha^*\rVert_\infty=m$. Interchanging the roles of $a$ and $b$ ($p=m$, $q=1$) gives the second vector, and the negatives are covered by the last statement of \cref{thm:posettwoisolates}. Conversely, let $\alpha$ be a primitive facet normal with $\alpha_a\alpha_b>0$ and $\lVert\alpha\rVert_\infty=m$. Replacing $\alpha$ by $-\alpha$ if necessary, \cref{thm:posettwoisolates} gives $\alpha=-\mathbf{e}_S+p\,\mathbf{e}_a+q\,\mathbf{e}_b$ with $p+q-1\leq|S|\leq m$ and $\max\{p,q\}=\lVert\alpha\rVert_\infty=m$. Because the smaller of $p,q$ is at least $1$, we get $p+q\geq m+1$, so $|S|=m=p+q-1$; that is, $S=C_m$ and $\{p,q\}=\{1,m\}$.
\end{proof}

For $P=P^{(m)}$, by \cref{cor:bijection}, $\alpha^*$ gives a vertical facet of $\D$ whose describing inequality \eqref{eq:FP} has $x_b$-coefficient $m=d-2$; as $\alpha^*$ is not a multiple of any $\pm\mathbf{e}_i$ or $\mathbf{e}_i-\mathbf{e}_j$\,, this facet is not described by the full optimal big-M lifting of a facet-describing inequality of any simplex of Stanley's triangulation. Thus, already for $\ell(P)=2$, the $x$-coefficients of the inequalities describing vertical facets can grow linearly in $d$, and by \cref{thm:sharpchaintwoelements} they grow at most linearly.

\begin{theorem}\label{thm:posetopposite}
Let $\alpha\in\mathbb{Z}^d$ with $\alpha_a\alpha_b<0$. Then $\alpha$ is a primitive facet normal of $M(P^{(m)})$ if and only if $\alpha_a=-\alpha_b\in\{-1,1\}$, every chain coefficient $\alpha_c$\,, $c\in C_m$\,, lies in $\{-1,0,1\}$, and the nonzero chain coefficients alternate in sign along $c_1\prec\cdots\prec c_m$; equivalently, $\alpha_a=-\alpha_b\in\{-1,1\}$ and the chain prefix sums $\sum_{j\leq i}\alpha_{c_j}$\,, $0\leq i\leq m$, take at most two values, which are consecutive integers. In particular, every such facet normal lies in $\{-1,0,1\}^d$, and there are exactly $2(2^{m+1}-1)$ of them.
\end{theorem}

\begin{proof}
    A sequence in $\{-1,0,1\}$ has alternating nonzero entries exactly when its prefix sums take at most two consecutive values, so the two forms of the condition agree. A chain part with $k\geq1$ nonzero entries is determined by its support and the sign of its first nonzero entry, which gives $1+\sum_{k\geq1}2\binom mk=2^{m+1}-1$ chain parts; the sign of $\alpha_a$ doubles this.

    \emph{Necessity: setup.} Interchanging $a$ and $b$ is an automorphism of $P^{(m)}$, so we may assume $\mu:=\alpha_a>0$ and $\nu:=-\alpha_b>0$. Let $R$ be the chain reversal from the proof of \cref{thm:posettwoisolates}, and let $\varsigma$ interchange the coordinates $a$ and $b$. Both $x\mapsto\mathbf{e}-Rx$ and $\varsigma$ permute the simplices of $\mathcal{T}(P^{(m)})$, so $\tilde\alpha:=-\varsigma R\alpha$ is again a primitive facet normal, with $\tilde\alpha_a=\nu$ and $\tilde\alpha_b=-\mu$; that is, $\tilde\alpha$ has the roles of $\mu$ and $\nu$ interchanged. By \cref{lem:dimofFMP}, the vectors orthogonal to all maximizing-block indicators are exactly the multiples of $\alpha$. Hence every maximizing block has $\alpha$-sum zero, and $(\dagger)$ no nonzero vector with vanishing $a$-coordinate is orthogonal to all block indicators.
    We use $\Sigma_i$\,, $\Sigma_{\min}$\,, $W_\pm$\,, and endpoint indices as in \cref{subsec:generalbound}, with $C=C_m$ and $E=\{a,b\}$. At a position of a linear extension where exactly $i$ chain elements have appeared, the prefix sum $\Phi^k_r$ equals $\Sigma_i$\,, plus $\mu$ if $a$ has appeared, minus $\nu$ if $b$ has appeared. For a maximizing block $B$ with endpoint indices $r\leq s$, let $\varepsilon_a\,,\varepsilon_b\in\{0,1\}$ indicate whether $a\in B$ and $b\in B$. Then
    \begin{equation}\label{eq:sruv}
    \Sigma_s-\Sigma_r+\varepsilon_a\mu-\varepsilon_b\nu=0.
    \end{equation}

    \emph{Step 1: endpoints.} Let $\Sigma_\lambda=\Sigma_{\min}$\,. Because $W_+(\alpha_E)=\mu$ and $W_-(\alpha_E)=\nu$, \cref{lem:window} gives, for every endpoint index $i$ of a maximizing block,
    \begin{equation}\label{eq:st}
    i\leq\lambda\ \Longrightarrow\ \Sigma_i\leq\Sigma_{\min}+\mu,
    \qquad
    i\geq\lambda\ \Longrightarrow\ \Sigma_i\leq\Sigma_{\min}+\nu.
    \end{equation}

    \emph{Step 2: $\nu\leq \mu$.} Suppose that $\nu>\mu$, and let $\lambda$ be the largest index with $\Sigma_\lambda=\Sigma_{\min}$\,. We show that a maximizing block $B$ contains $b$ if and only if $r\leq\lambda<s$. If $b\in B$ and $a\notin B$, then $\Sigma_s=\Sigma_r+\nu$; if $s\leq\lambda$, \eqref{eq:st} would give $\Sigma_r\leq\Sigma_{\min}+\mu-\nu<\Sigma_{\min}$\,, so $s>\lambda$, and then \eqref{eq:st} gives $\Sigma_r\leq\Sigma_{\min}$\,, so $\Sigma_r=\Sigma_{\min}$ and $r\leq\lambda$. If $a,b\in B$, then only chain elements lie before and after $B$, so $\Sigma_r=\min_{i\leq r}\Sigma_i$ and $\Sigma_s=\min_{i\geq s}\Sigma_i$\,, while $\Sigma_s=\Sigma_r+\nu-\mu>\Sigma_r$; now $r>\lambda$ would give $\Sigma_r=\Sigma_{\min}$\,, contradicting the choice of $\lambda$, and $s\leq\lambda$ would give $\Sigma_s=\Sigma_{\min}\leq\Sigma_r$\,. Finally, if $b\notin B$ and $r\leq\lambda<s$, then $\Sigma_s>\Sigma_{\min}$ by the choice of $\lambda$, so $\Sigma_r=\Sigma_s+\varepsilon_a\mu>\Sigma_{\min}$ and $r<\lambda$. The position right after $c_\lambda$ then lies strictly inside $B$, and its prefix sum exceeds that at the first endpoint of $B$ by at most $(\Sigma_{\min}-\Sigma_r)+\varepsilon_a\mu=\Sigma_{\min}-\Sigma_s<0$, contradicting minimality. This proves the claim. Consequently $\mathbf{e}_b-\mathbf{e}_{c_{\lambda+1}}$ (or $\mathbf{e}_b$ if $\lambda=m$) is orthogonal to all block indicators, contradicting $(\dagger)$. Hence $\nu\leq \mu$, and applying this to $\tilde\alpha$ gives $\mu\leq \nu$, so $\mu=\nu$.

    \emph{Step 3: the chain coefficients.} Because $\mu=\nu$, \eqref{eq:st} gives $\Sigma_r\,,\Sigma_s\in[\Sigma_{\min}\,,\Sigma_{\min}+\mu]$ for every maximizing block, and by \eqref{eq:sruv} either $\Sigma_r=\Sigma_s$\,, or $|\Sigma_s-\Sigma_r|=\mu$ and $\{\Sigma_r\,,\Sigma_s\}=\{\Sigma_{\min}\,,\Sigma_{\min}+\mu\}$. Let $\eta_i:=1$ if $\Sigma_i\in\{\Sigma_{\min}\,,\Sigma_{\min}+\mu\}$ and $\eta_i:=0$ otherwise; then $\eta_r=\eta_s$ for every maximizing block. The vector $\gamma$ with $\gamma_{c_i}:=\eta_i-\eta_{i-1}$ and $\gamma_a=\gamma_b:=0$ satisfies $\gamma^\top\mathbf{e}_B=\eta_s-\eta_r=0$, so $\gamma=0$ by $(\dagger)$. Thus $\eta$ is constant, equal to $\eta_\lambda=1$, and every $\Sigma_i$ lies in $\{\Sigma_{\min}\,,\Sigma_{\min}+\mu\}$. Hence the chain coefficients lie in $\{-\mu,0,\mu\}$ and the nonzero ones alternate in sign; as $\alpha$ is primitive, $\mu=1$.

    \emph{Sufficiency.} We may assume that $\alpha_a=1$ and $\alpha_b=-1$, so that $\alpha$ is primitive and every $\Sigma_i$ lies in $\{\Sigma_{\min}\,,\Sigma_{\min}+1\}$. Each linear extension below inserts $a$ and $b$ next to a single chain element; its prefix sums are chain prefix sums except at most at two positions next to the inserted elements, and one checks in each case that these are at least $\Sigma_{\min}$ as well. So the minimum is $\Sigma_{\min}$\,, and the positions listed with value $\Sigma_{\min}$ are minimal.
    \begin{enumerate}[label=(\roman*)]
        \item Insert $a,b$ right after some $c_{i_0}$ with $\Sigma_{i_0}=\Sigma_{\min}$ (at the beginning if $i_0=0$). The values $\Sigma_{\min}\,,\Sigma_{\min}+1,\Sigma_{\min}$ make $\{a,b\}$ a maximizing block.
        \item If $\alpha_{c_j}=1$, insert $b,a$ right after $c_j$; the values $\Sigma_{\min}$ before $c_j$\,, $\Sigma_{\min}+1$ after it, and $\Sigma_{\min}$ after $b$ make $\{c_j\,,b\}$ a block.
        \item If $\alpha_{c_j}=-1$, insert $b,a$ right before $c_j$; the values $\Sigma_{\min}$ after $b$, $\Sigma_{\min}+1$ after $a$, and $\Sigma_{\min}$ after $c_j$ make $\{a,c_j\}$ a block.
        \item If $\alpha_{c_j}=0$, insert $a,b$ right after $c_j$ when $\Sigma_j=\Sigma_{\min}$\,, and $b$ right before and $a$ right after $c_j$ when $\Sigma_j=\Sigma_{\min}+1$; either way the prefix sums just before and after $c_j$ equal $\Sigma_{\min}$\,, so $\{c_j\}$ is a block.
    \end{enumerate}
    We obtain $\mathbf{e}_a+\mathbf{e}_b$ and, for each $j\in[m]$, a block indicator whose support contains $c_j$ and lies in $\{c_j\,,a,b\}$. These $m+1=d-1$ vectors are linearly independent, because each of the latter has a chain coordinate that vanishes in all the others, so $\alpha$ is a facet normal by \cref{lem:dimofFMP}.
\end{proof}

\begin{proposition}\label{lem:twoisolateszero}
    A vector $\alpha\in\mathbb{Z}^d$ with $\alpha_a\alpha_b=0$ is a primitive facet normal of $M(P^{(m)})$ if and only if it is one of
    \[
    \begin{gathered}
    -\mathbf{e}_{c_1}\,,\quad \mathbf{e}_{c_m}\,,\quad \pm\mathbf{e}_a\,,\quad \pm\mathbf{e}_b\,,\quad \mathbf{e}_{c_i}-\mathbf{e}_{c_{i+1}}\ (i\in[m-1]),\\
    \pm(\mathbf{e}_a-\mathbf{e}_c),\quad \pm(\mathbf{e}_b-\mathbf{e}_c)\quad (c\in C_m).
    \end{gathered}
    \]
    In particular, there are exactly $5(m+1)$ of them, and each has $\lVert\alpha\rVert_\infty=1$.
\end{proposition}

\begin{proof}
    Each listed vector is a facet normal of some simplex of $\mathcal{T}(P^{(m)})$, because $\Delta_\sigma$ has the facet normals $-\mathbf{e}_{\sigma_1}$\,, $\mathbf{e}_{\sigma_d}$ and $\mathbf{e}_{\sigma_t}-\mathbf{e}_{\sigma_{t+1}}$ ($t\in[d-1]$): take a linear extension $\sigma$ that starts with $c_1$\,, $a$ or $b$; ends with $c_m$\,, $a$ or $b$; has $c_i\,,c_{i+1}$ adjacent; or has $a$ (or $b$) immediately before or immediately after $c$. By \cref{lem:persistenceaddition}, each listed vector is therefore a facet normal of $M(P^{(m)})$, and it is primitive.

    Conversely, let $\alpha$ be a primitive facet normal with $\alpha_a\alpha_b=0$. Interchanging $a$ and $b$ if necessary, we may assume that $\alpha_b=0$. By the second statement of \cref{lem:deletezeros} with $Z=\{b\}$, the restriction $\overline\alpha$ of $\alpha$ to $C_m\cup\{a\}$ is a primitive facet normal of $M(\overline P)$, where $\overline P:=P^{(m)}|_{C_m\cup\{a\}}$\,. Deleting $a$ from $\overline P$ leaves a chain, so by \cref{thm:verticalarelifting} every facet normal of $M(\overline P)$ is a facet normal of some simplex $\Delta_\rho$ of $\mathcal{T}(\overline P)$; because $\overline\alpha$ is primitive, it is $-\mathbf{e}_{\rho_1}$\,, $\mathbf{e}_{\rho_{m+1}}$\,, or $\mathbf{e}_{\rho_t}-\mathbf{e}_{\rho_{t+1}}$\,. In a linear extension $\rho$ of $\overline P$, the first element is $c_1$ or $a$, the last is $c_m$ or $a$, and two consecutive elements are $c_i\,,c_{i+1}$ or involve $a$. Appending $\alpha_b=0$ shows that $\alpha$ is in the list. Finally, the list has $2+4+(m-1)+2m+2m=5(m+1)$ entries.
\end{proof}

\Cref{thm:posettwoisolates,thm:posetopposite,lem:twoisolateszero} together give a complete classification. In the count below, $\varphi$ denotes Euler's totient function.

\begin{corollary}\label{cor:twoisolatescount}
    Let $m\geq1$. The primitive facet normals of $M(P^{(m)})$ are exactly the vectors $\pm(-\mathbf{e}_S+p\,\mathbf{e}_a+q\,\mathbf{e}_b)$ of \cref{thm:posettwoisolates}, the vectors of \cref{thm:posetopposite}, and the vectors of \cref{lem:twoisolateszero}. Their number, which by \cref{cor:bijection} is the number of vertical facets of $\D$ for $P=P^{(m)}$, equals
    \[
    2\textstyle\sum_{s=2}^{m+1}\varphi(s)\bigl(\tbinom{m}{s-1}+\tbinom{m}{s}\bigr)+2^{m+2}-2+5(m+1).
    \]
\end{corollary}

\begin{proof}
    The three results cover the cases $\alpha_a\alpha_b>0$, $\alpha_a\alpha_b<0$ and $\alpha_a\alpha_b=0$. For the count, in \cref{thm:posettwoisolates} put $s:=p+q$; the coprime pairs $(p,q)$ of positive integers with $p+q=s\geq2$ number $\varphi(s)$, because $\gcd(p,s-p)=\gcd(p,s)$, and for each of them $S$ is any subset of $C_m$ of size $s-1$ or $s$, which requires $s\leq m+1$; the sign doubles this. \Cref{thm:posetopposite} contributes $2(2^{m+1}-1)$ and \cref{lem:twoisolateszero} contributes $5(m+1)$ normals.
\end{proof}

Among these normals, those of alcoved type are the $5(m+1)$ vectors of \cref{lem:twoisolateszero} together with $\pm(\mathbf{e}_a-\mathbf{e}_b)$ from \cref{thm:posetopposite}, $5m+7$ in all; this agrees with \cref{cor:alcovedlist}, as $P^{(m)}$ has $\min(P^{(m)})=\max(P^{(m)})=3$, $\operatorname{inc}(P^{(m)})=2m+1$, and $\operatorname{cov}(P^{(m)})=m-1$. For $m=1$, the total count is $18$, matching the facet normals $\pm\mathbf{e}_i$\,, $\mathbf{e}_i-\mathbf{e}_j$\,, and $\pm(\mathbf{e}_i+\mathbf{e}_j-\mathbf{e}_k)$ of $M(P)$ for the antichain on three elements (\cref{ex:antichain}).

\begin{theorem}\label{thm:sharpchaintwoelements}
Let $d\geq3$, and let $P$ be a poset on $V=C\sqcup\{a,b\}$, where $C$ is a chain of $P$ with $d-2$ elements; equivalently, $\ell(P)\leq2$. Then every primitive facet normal $\alpha$ of $M(P)$ satisfies $\lVert\alpha\rVert_\infty\leq d-2$. The bound is attained by $P=P^{(d-2)}$; for this poset and $d\geq4$, equality holds exactly for the four vectors $\pm(-\mathbf{e}_{C}+\mathbf{e}_a+(d-2)\,\mathbf{e}_b)$ and $\pm(-\mathbf{e}_{C}+(d-2)\,\mathbf{e}_a+\mathbf{e}_b)$.
\end{theorem}

\begin{proof}
    Put $m:=d-2\geq1$ and label $C$ as $c_1\prec\cdots\prec c_m$\,. As noted at the beginning of \cref{subsec:chaintwo}, every facet normal of $M(P)$ is a facet normal of $M(P^{(m)})$; because the primitive integer vector on a rational ray is unique, it suffices to bound the primitive facet normals $\alpha$ of $M(P^{(m)})$. If $\alpha_a\alpha_b>0$, then $\lVert\alpha\rVert_\infty=\max\{p,q\}\leq m$ by \cref{thm:posettwoisolates} (see the paragraph following it). If $\alpha_a\alpha_b<0$, then $\alpha\in\{-1,0,1\}^d$ by \cref{thm:posetopposite}, and if $\alpha_a\alpha_b=0$, then $\lVert\alpha\rVert_\infty=1$ by \cref{lem:twoisolateszero}. Because $m\geq1$, in every case $\lVert\alpha\rVert_\infty\leq m=d-2$.

    For $P=P^{(m)}$, \cref{cor:posetsharp} gives normals with $\lVert\alpha\rVert_\infty=m$, so the bound is attained. If $m\geq2$, the normals with $\alpha_a\alpha_b\leq0$ have $\lVert\alpha\rVert_\infty=1<m$, so equality can only occur when $\alpha_a\alpha_b>0$, and \cref{cor:posetsharp} identifies those normals as the four listed vectors.
\end{proof}

\subsection{Lower bounds, and the growth for fixed \texorpdfstring{$\ell(P)$}{l(P)}}\label{subsec:lowerbound}
By \cref{subsec:chaintwo}, $H_\ell(d)$ is attained by $P^{(d-\ell,\ell)}$\,. \Cref{thm:posetheight} gives $H_\ell(d)\leq\delta_\ell\binom{d-1}{\ell-1}$\,, of degree $\ell-1$ in $d$, and \cref{thm:sharpchaintwoelements} gives $H_2(d)=d-2$. We show that the degree $\ell-1$ is attained for every fixed $\ell\geq2$, that the entries can be exponential in $d$ when $\ell$ grows with $d$, and we give a better quadratic lower bound for $\ell=3$. All constructions rest on the following lemma.

\begin{lemma}\label{lem:interval}
Let $P=P^{(m,\ell)}$, and let $\alpha\in\mathbb{R}^d$ satisfy $\alpha_c<0$ for every chain element $c$, $\alpha_e>0$ for every isolated element $e$, and $\sum_{e}\alpha_e>\sum_{c}|\alpha_c|$. If $I=\{c_{r+1}\,,\ldots,c_s\}$ is an interval of the chain and $A$ is a nonempty set of isolated elements with $\sum_{e\in A}\alpha_e=\sum_{c\in I}|\alpha_c|$, then $I\cup A$ is a maximizing block of some linear extension of $P$.
\end{lemma}

\begin{proof}
Because $\sum_{e\in A}\alpha_e\leq\sum_c|\alpha_c|<\sum_e\alpha_e$\,, $A$ is not all of the isolated elements, and $I\neq\varnothing$. Consider the linear extension that lists $c_1\,,\ldots,c_r$\,, then the elements of $A$, then $c_{r+1}\,,\ldots,c_s$\,, then the remaining isolated elements, and finally $c_{s+1}\,,\ldots,c_m$\,. With $\Sigma_r:=\sum_{i\leq r}\alpha_{c_i}\leq0$ as in \cref{subsec:generalbound}, its prefix sums decrease from $0$ to $\Sigma_r$\,, exceed $\Sigma_r$ strictly inside the run $A\cup I$ and return to $\Sigma_r$ at its end, then increase along the remaining isolated elements and decrease to $\sum_e\alpha_e-\sum_c|\alpha_c|>0\geq\Sigma_r$\,. So the two positions flanking $A\cup I$ are consecutive minimal positions.
\end{proof}

\begin{theorem}[Product construction]\label{thm:product}
Let $\ell\geq2$, let $k_1\,,\ldots,k_{\ell-1}\geq2$ be integers, and put $u_0:=1$, $u_i:=k_iu_{i-1}$ for $1\leq i\leq\ell-1$, and $d:=2+\sum_{i=1}^{\ell-1}k_i$\,. Let $P=P^{(d-\ell,\ell)}$ have isolated elements $a_0\,,\ldots,a_{\ell-1}$\,, and partition its chain, in chain order, into consecutive groups $G_1\,,\ldots,G_{\ell-1}$ with $|G_1|=k_1$ and $|G_i|=k_i-1$ for $i\geq2$. Define $\alpha$ by $\alpha_{a_i}:=u_i$ and $\alpha_c:=-u_{i-1}$ for $c\in G_i$\,. Then $\alpha$ is a primitive facet normal of $M(P)$, and $\lVert\alpha\rVert_\infty=u_{\ell-1}=\prod_{i=1}^{\ell-1}k_i$\,. \end{theorem}

\begin{proof}
The chain has $k_1+\sum_{i\geq2}(k_i-1)=d-\ell\geq\ell$ elements, so $\ell(P)=\ell$. Let $C_i:=G_1\cup\cdots\cup G_i$\,. By induction, $\sum_{c\in C_i}|\alpha_c|=u_i$\,: this is $k_1u_0=u_1$ for $i=1$, and $u_{i-1}+(k_i-1)u_{i-1}=k_iu_{i-1}=u_i$ in general. So the total chain weight $u_{\ell-1}$ is less than $\sum_iu_i$\,, and by \cref{lem:interval} the following $d-1$ sets are maximizing blocks: $\{c,a_{i-1}\}$ for every $c\in G_i$ and every $i$, and $C_i\cup\{a_i\}$ for $1\leq i\leq\ell-1$.

Let $y$ be orthogonal to their indicator vectors. The blocks $\{c,a_{i-1}\}$ give $y_c=-y_{a_{i-1}}$ for $c\in G_i$\,. The block $C_1\cup\{a_1\}$ gives $y_{a_1}=k_1y_{a_0}$\,, and for $i\geq2$, the block $C_i\cup\{a_i\}$ gives $y_{a_i}=-\sum_{c\in C_i}y_c=y_{a_{i-1}}+(k_i-1)y_{a_{i-1}}=k_iy_{a_{i-1}}$\,. Hence $y=y_{a_0}\alpha$, so the indicator vectors span $\alpha^\perp$, and $\alpha$ is a facet normal by \cref{lem:dimofFMP}. It is primitive because $\alpha_{a_0}=1$, and its largest entry is $u_{\ell-1}$ because the chain entries have absolute value at most $u_{\ell-2}<u_{\ell-1}$\,.
\end{proof}

\begin{corollary}\label{cor:growth}
Let $\ell\geq2$ and $d\geq2\ell$, and write $d-2=q(\ell-1)+r$ with $0\leq r<\ell-1$. Then $q^{\ell-1-r}(q+1)^r\leq H_\ell(d)\leq\delta_\ell\binom{d-1}{\ell-1}$\,. Consequently, for every fixed $\ell\geq2$, $H_\ell(d)=\Theta(d^{\ell-1})$, with
\[
(\ell-1)^{-(\ell-1)}\;\leq\;\liminf_{d\to\infty}\frac{H_\ell(d)}{d^{\ell-1}}\;\leq\;\limsup_{d\to\infty}\frac{H_\ell(d)}{d^{\ell-1}}\;\leq\;\frac{\delta_\ell}{(\ell-1)!}\;\leq\;\frac{(\ell+1)^{(\ell+1)/2}}{2^\ell\,(\ell-1)!}\,.
\]
\end{corollary}

\begin{proof}
Because $d-2\geq2(\ell-1)$, we have $q\geq2$. Apply \cref{thm:product} with $\ell-1-r$ factors equal to $q$ and $r$ factors equal to $q+1$; their sum is $d-2$ and their product is $q^{\ell-1-r}(q+1)^r$. The upper bound and the first bound on the $\limsup$ are \cref{thm:posetheight}, and the last inequality is \cref{lem:delta4}. Finally, $q\geq(d-\ell)/(\ell-1)$, so $q^{\ell-1}/d^{\ell-1}\to(\ell-1)^{-(\ell-1)}$.
\end{proof}

For $\ell=2$, \cref{thm:product} gives $d-2$, which is sharp by \cref{thm:sharpchaintwoelements}. Letting $\ell$ grow with $d$ gives exponential entries.

\begin{corollary}\label{cor:exp}
For every $d\geq6$, put $s:=d-2$ and
\[
\mathfrak{B}(s):=\begin{cases}3^{s/3},&s\equiv0\pmod3,\\ 4\cdot3^{(s-4)/3},&s\equiv1\pmod3,\\ 2\cdot3^{(s-2)/3},&s\equiv2\pmod3.\end{cases}
\]
For $\ell:=\lfloor(d-1)/3\rfloor+1$, the poset $P=P^{(d-\ell,\ell)}$ has a primitive facet normal $\alpha$ of $M(P)$ with $\lVert\alpha\rVert_\infty=\mathfrak{B}(d-2)$, and hence, by \cref{cor:bijection}, $\D$ has a vertical facet whose describing inequality has an $x$-coefficient of this size; $\mathfrak{B}(d-2)=\Theta(3^{d/3})$.
\end{corollary}

\begin{proof}
Apply \cref{thm:product} with all factors equal to $3$ if $s\equiv0$; with one factor $4$ and the others $3$ if $s\equiv1$; and with one factor $2$ and the others $3$ if $s\equiv2 \pmod 3$. The factors sum to $s$ and multiply to $\mathfrak{B}(s)$; their number is $\lfloor(s+1)/3\rfloor=\ell-1$.
\end{proof}

The number $\mathfrak{B}(s)$ is the largest product of positive integers with sum at most $s$, and it has the same meaning and notation in \cite{QL_staircase}: there, $\mathfrak{B}(d-2)$ is the upper bound for the entries of a primitive facet normal of a Minkowski sum of at most three staircase simplices \SF{\SFsize(b)}. \Cref{cor:exp} attains that bound exactly, but with a sum of $|\mathcal L(P)|$ simplices, which is far more than three.

By contrast, for fixed $\ell$ the growth is polynomial (\cref{cor:growth}), and the exponential entries of \cite{QL_staircase} for sums of arbitrary staircase simplices occur already with three summands. For $\ell=3$ and $d\geq6$, \cref{thm:product} gives $H_3(d)\geq\lfloor(d-2)^2/4\rfloor$; the following construction does better.

\begin{theorem}[Coprime-grid construction]\label{thm:grid}
Let $p,q,L\geq1$ be integers with $\gcd(p,q)=1$, and put $d:=p+q+L+2$ and $K:=pq+L(p+q)$. Let $0=t_0<t_1<\cdots<t_{p+q-1}=pq$ be the elements of $\{0,q,2q,\ldots,pq\}\cup\{0,p,2p,\ldots,pq\}$, and let $P=P^{(d-3,3)}$ have isolated elements $u,v,w$. Define $\alpha$ by $\alpha_{c_j}:=-(t_j-t_{j-1})$ for $1\leq j\leq p+q-1$, $\alpha_{c_j}:=-(p+q)$ for $p+q\leq j\leq d-3$, and $\alpha_u:=K$, $\alpha_v:=q$, $\alpha_w:=p$. Then $\alpha$ is a primitive facet normal of $M(P)$ with $\lVert\alpha\rVert_\infty=K$.
\end{theorem}

\begin{proof}
Because $\gcd(p,q)=1$, the two arithmetic progressions meet only in $0$ and $pq$, so there are indeed $p+q-1$ gaps, and the chain has $p+q-1+L=d-3$ elements. The chain entries are negative with total $-(pq+L(p+q))=-K$, and the isolated entries are positive with sum $K+p+q>K$, so \cref{lem:interval} applies. It shows that the following $p+q+L+1=d-1$ sets are maximizing blocks: for each of the $p$ intervals $[jq,(j+1)q]$ of the first progression, the chain elements whose gaps lie in it, together with $v$; for each of the $q$ intervals $[ip,(i+1)p]$ of the second progression, the corresponding chain elements together with $w$; each of the $L$ chain elements $c_j$ with $j\geq p+q$, together with $v$ and $w$; and the whole chain together with $u$.

Let $y$ be orthogonal to their indicator vectors, and let $Y(t_j):=\sum_{i\leq j}y_{c_i}$\,, with $Y(0)=0$. Summing the $p$ equations of the first kind gives $Y(pq)=-py_v$\,, and the $q$ equations of the second kind give $Y(pq)=-qy_w$\,; hence $y_v=qt$ and $y_w=pt$ for $t:=y_w/p$. The equations then give $Y(jq)=-jqt$ and $Y(ip)=-ipt$, that is, $Y(t_j)=-t_jt$ for every $j$, because every $t_j$ lies in one of the two progressions; so $y_{c_j}=-(t_j-t_{j-1})t$ for $j\leq p+q-1$. The equations of the third kind give $y_{c_j}=-(p+q)t$ for $j\geq p+q$, and the last one gives $y_u=Kt$. Thus $y=t\alpha$, the $d-1$ indicator vectors span $\alpha^\perp$, and $\alpha$ is a facet normal by \cref{lem:dimofFMP}. It is primitive because $\gcd(\alpha_v\,,\alpha_w)=\gcd(p,q)=1$, and its largest entry is $K$, as $K\geq p+q$ and every gap is at most $\max(p,q)$.
\end{proof}

\begin{corollary}\label{cor:grid}
Let $d\geq5$ and $s:=d-2$. Then
\[
H_3(d)\;\geq\;\begin{cases}s^2/3-1,&s\equiv0\pmod3\text{ and }s\geq6,\\ \lfloor s^2/3\rfloor,&\text{otherwise}.\end{cases}
\]
In particular $H_3(d)\geq\lfloor(d-2)^2/3\rfloor-1$, and $\liminf_{d\to\infty}H_3(d)/d^2\geq1/3$.
\end{corollary}

\begin{proof}
In \cref{thm:grid}, $s=p+q+L$. For $s=3k+1$ take $(p,q,L)=(k,k+1,k)$, giving $K=3k^2+2k=\lfloor s^2/3\rfloor$; for $s=3k+2$ take $(k,k+1,k+1)$, giving $K=3k^2+4k+1=\lfloor s^2/3\rfloor$; for $s=3k$ with $k\geq2$ take $(k-1,k,k+1)$, giving $K=3k^2-1=s^2/3-1$; and for $s=3$ take $(1,1,1)$, giving $K=3=\lfloor s^2/3\rfloor$.
\end{proof}

\begin{remark}\label[remark]{rem:sixteen}
For $d=5,6,7,8$, \cref{cor:grid} gives $3,5,8,11$, which are the largest entries found by computation (\cref{rem:open}); for $d=5,6$ our computations enumerate all facet normals (complete by \cref{thm:upper3}), and they give $H_3(5)=3$ and $H_3(6)=5$. For $d=9$ it gives $16$: with $(p,q,L)=(2,3,2)$, the progressions $\{0,3,6\}$ and $\{0,2,4,6\}$ have gaps $2,1,1,2$, and
\[
\alpha=-2\,\mathbf{e}_{c_1}-\mathbf{e}_{c_2}-\mathbf{e}_{c_3}-2\,\mathbf{e}_{c_4}-5\,\mathbf{e}_{c_5}-5\,\mathbf{e}_{c_6}+16\,\mathbf{e}_u+3\,\mathbf{e}_v+2\,\mathbf{e}_w
\]
is a primitive facet normal of $M(P^{(6,3)})$; by \cref{cor:bijection}, the corresponding vertical facet of $\D$ has $x_u$-coefficient $16$.
\end{remark}

For $\ell=4$, the following vector gives $H_4(d)\geq4d-19$, which for $7\leq d\leq9$ exceeds what \cref{thm:product} and, via \cref{lem:persistenceaddition}, \cref{cor:grid} provide.

\begin{proposition}\label{prop:lowerboundfour}
    Let $m\geq2$, $d:=m+4$, and let $P=P^{(m,4)}$ have chain $c_1\prec\cdots\prec c_m$ and isolated elements $a_1\,,a_2\,,b,e$. Then
    \[
    \alpha:=-4\textstyle\sum_{i=1}^{m-1}\mathbf{e}_{c_i}-3\,\mathbf{e}_{c_m}+\mathbf{e}_{a_1}+\mathbf{e}_{a_2}+2\,\mathbf{e}_b+(4m-3)\,\mathbf{e}_e
    \]
    is a primitive facet normal of $M(P)$ with $\lVert\alpha\rVert_\infty=4m-3=4d-19$.
\end{proposition}

\begin{proof}
    The chain entries are negative with total $-(4m-1)$, and the isolated entries are positive with sum $4m+1$, so by \cref{lem:interval} the following $d-1$ sets are maximizing blocks: $\{a_1\,,a_2\,,b,c_i\}$ for $i\in[m-1]$, $\{a_j\,,b,c_m\}$ for $j\in\{1,2\}$, $\{e,b,c_1\,,\ldots,c_m\}$, and $\{e,a_1\,,a_2\,,c_1\,,\ldots,c_m\}$. If $y$ is orthogonal to their indicator vectors, then $y_{c_i}=-S$ for $i<m$ with $S:=y_{a_1}+y_{a_2}+y_b$\,, $y_{a_1}=y_{a_2}=:t$ and $y_{c_m}=-t-y_b$\,, and the last two blocks give $y_b=2t$; hence $y=t\alpha$, and $\alpha$ is a primitive facet normal with largest entry $4m-3$ by \cref{lem:dimofFMP}.
\end{proof}

\subsection{Sharper constants for \texorpdfstring{$\ell(P)=3$ and $\ell(P)=4$}{l(P)=3 and l(P)=4}}\label{subsec:upper3}
For $\ell=3$, \cref{thm:posetheight} gives $H_3(d)\leq(d-1)(d-2)$, because $\delta_3=2$; we improve the constant.

\begin{theorem}\label{thm:upper3}
Let $P$ be a poset on $[d]$ with $\ell(P)=3$. Every primitive facet normal $\alpha$ of $M(P)$ satisfies $\lVert\alpha\rVert_\infty\leq\lfloor(d-1)^2/2\rfloor$.
\end{theorem}

\begin{proof}
Let $C$, $E$, and $m$ be as in \cref{subsec:generalbound}, with $C$ a longest chain, so that $m=d-3$ and $|E|=3$. We may assume $\lVert\alpha\rVert_\infty\geq2$. Then \cref{lem:window} gives $|\alpha_c|\leq W(\alpha_E)$ for every $c\in C$; in particular $\alpha_E\neq0$. Also $|\alpha_e|\leq W(\alpha_E)\in\mathbb{Z}$ for $e\in E$, so it suffices to show $W(\alpha_E)\leq(d-1)^2/2$. We order the coordinates so that the elements of $E$ come last; this changes the cofactor vector of \cref{lem:crossproduct} at most by a global sign.

By \cref{lem:dimofFMP}, choose $d-1$ maximizing blocks $B_1\,,\ldots,B_{d-1}$ of linear extensions of $P$ whose indicator vectors are linearly independent, and let $A$ be the $(d-1)\times d$ matrix with rows $\mathbf{e}_{B_j}^\top$\,; then $\ker A=\Span\{\gamma\}$ for the cofactor vector $\gamma$ of \cref{lem:crossproduct}, and $\alpha=\pm\gamma/g$ with $g\geq1$. As $\alpha_E\neq0$, the $m$ columns of $A$ indexed by $C$ are linearly independent, and $B_j\cap C$ is an interval $\{c_{r_j+1}\,,\ldots,c_{s_j}\}$ with $0\leq r_j\leq s_j\leq m$ (empty if $r_j=s_j$). Write $p_j:=(\mathbf{e}_{B_j})_E\in\{0,1\}^E$ and let $\Pi$ be the $(d-1)\times3$ matrix with rows $p_j^\top$\,.

Replace the columns $A_{c_1}\,,\ldots,A_{c_m}$ of $A$ by $A_{c_1}-A_{c_2}\,,\ldots,A_{c_{m-1}}-A_{c_m}\,,A_{c_m}$\,, a unimodular operation on the chain columns that leaves the $E$ columns unchanged. Row $j$ of the new chain part has entry $1$ in column $s_j$ if $s_j\geq1$, entry $-1$ in column $r_j$ if $r_j\geq1$, and zeros elsewhere; so the new chain part is the reduced incidence matrix $\Lambda$ (vertex $0$ deleted) of the directed multigraph $\Gamma$ on the vertex set $\{0,1,\ldots,m\}$ with the $d-1$ arcs $r_j\to s_j$\,, loops allowed. The columns of $\Lambda$ are linearly independent, so $\Gamma$ is connected; as $\Gamma$ has $d-1=m+2$ arcs on $m+1$ vertices, its cycle rank is $2$. Choose a spanning tree $\Upsilon$ of $\Gamma$; the $m$ rows of $\Lambda$ indexed by the arcs of $\Upsilon$ form a matrix $\Lambda_\Upsilon$ with $\det\Lambda_\Upsilon=\pm1$ (expand along the column of a leaf of $\Upsilon$ other than $0$, which contains a single nonzero entry $\pm1$, and use induction). Let $f_1\,,f_2$ be the two arcs outside $\Upsilon$, and let $\nu_1\,,\nu_2\in\{-1,0,1\}^{d-1}$ be the signed indicator vectors of their fundamental cycles $Z_1\,,Z_2$ with respect to $\Upsilon$, so that $\nu_i^\top\Lambda=0$, $\nu_i(f_i)=1$, and $\nu_i(f_{3-i})=0$. Put $\xi:=\Pi^\top\nu_1$ and $\zeta:=\Pi^\top\nu_2$\,, the signed sums of the vectors $p_j$ along $Z_1$ and $Z_2$\,.

For $e\in E$, let $A^{(e)}$ be $A$ with column $e$ deleted; the column operation does not change $|\det A^{(e)}|=|\gamma_e|$. Order the rows of the transformed $A^{(e)}$ with those of $\Upsilon$ first; subtracting from the two remaining rows the combination $\Lambda_R\Lambda_\Upsilon^{-1}$ of the tree rows (where $\Lambda_R$ consists of the rows $f_1\,,f_2$ of $\Lambda$) annihilates their chain part, and the rows of $\bigl(-\Lambda_R\Lambda_\Upsilon^{-1}\ \ I_2\bigr)$ are exactly $\nu_1^\top,\nu_2^\top$, because the left kernel of $\Lambda$ is two-dimensional and $\nu_1\,,\nu_2$ are its unique elements with these entries on $f_1\,,f_2$\,. Hence
\[
|\gamma_e|=\bigl|\det\Lambda_\Upsilon\bigr|\cdot\bigl|\det\bigl(\nu_i^\top\Pi_{E\setminus\{e\}}\bigr)_{i=1,2}\bigr|=|(\xi\times\zeta)_e|,
\]
and because the signs are consistent across the three coordinates, $\gamma_E=\pm\,\xi\times\zeta$. As $\gamma\neq0$ is parallel to $\alpha$ and $\alpha_E\neq0$, we get $\xi\times\zeta\neq0$ and $W(\alpha_E)=W(\gamma_E)/g\leq W(\xi\times\zeta)$.

For every nonempty $S\subset E$ there is a linear map $T_S:\mathbb{R}^E\to\mathbb{R}^2$ with
\begin{equation}\label{eq:TS}
\textstyle\sum_{e\in S}(x\times y)_e=\det(T_Sx,T_Sy)\quad(x,y\in\mathbb{R}^E)\quad\text{and}\quad T_S(\{0,1\}^E)\subset[-1,1]^2 :
\end{equation}
writing $E=\{1,2,3\}$, take $T_{\{1\}}(x)=(x_2\,,x_3)$, $T_{\{1,2\}}(x)=(x_2-x_1\,,x_3)$, $T_{\{1,2,3\}}(x)=(x_1-x_3\,,x_2-x_3)$, and the cyclic permutations of these for the remaining sets. Because $\xi$ is a signed sum of $|Z_1|$ vectors in $\{0,1\}^E$, $T_S\xi\in[-|Z_1|,|Z_1|]^2$, and similarly for $\zeta$ and for signed sums along any set of arcs.

We now determine the shape of $Z_1$ and $Z_2$\,. Delete from $\Gamma$ its bridges and the vertices thereby isolated; this changes neither the cycle rank nor the fundamental cycles (every bridge belongs to $\Upsilon$ and to no cycle), the remaining graph is bridgeless, and every remaining vertex has degree at least two. Each remaining component is bridgeless with at least one arc and so has cycle rank at least one; as the ranks add up to $2$, at most two components remain. If two remain, each has cycle rank one, hence as many arcs as vertices, so the sum of $\deg v-2$ over its vertices is $0$; all its degrees equal two, and it is a cycle (a loop or a pair of parallel arcs included). Then $\Upsilon$ omits one arc from each, and $Z_1\,,Z_2$ are these two arc-disjoint cycles. If one component $\Gamma'$ remains, it has one more arc than vertices, so the sum of $\deg v-2$ over its vertices is $2$. Either one vertex has degree four and all others degree two, in which case $\Gamma'$ consists of two cycles sharing that vertex, and again $Z_1\,,Z_2$ are these arc-disjoint cycles; or two vertices $o_1\,,o_2$ have degree three and all others degree two. In the latter case every maximal path leaving $o_1$ through vertices of degree two ends at $o_1$ or $o_2$\,, and it cannot return to $o_1$\,, because then the third arc at $o_1$ would start a path to $o_2$ forming a bridge; hence $\Gamma'$ is a theta graph, the union of three internally disjoint paths $\pi_1\,,\pi_2\,,\pi_3$ from $o_1$ to $o_2$\,. The two arcs $f_1\,,f_2$ outside $\Upsilon$ then lie on different paths, as $\Upsilon$ contains no cycle; say $f_1\in\pi_1$ and $f_2\in\pi_3$\,, so that $Z_1=\pi_1\cup\pi_2$ and $Z_2=\pi_2\cup\pi_3$\,.

In the two-cycle case, with $q_1:=|Z_1|$ and $q_2:=|Z_2|$, \eqref{eq:TS} gives, for every $S$,
\[
\Bigl|\textstyle\sum_{e\in S}(\xi\times\zeta)_e\Bigr|=|\det(T_S\xi,T_S\zeta)|\leq2q_1q_2\leq(q_1+q_2)^2/2\leq(d-1)^2/2.
\]
In the theta case, orient the three paths from $o_1$ to $o_2$\,, let $\xi_1\,,\xi_2\,,\xi_3$ be the signed sums of the $p_j$ along $\pi_1\,,\pi_2\,,\pi_3$\,, and let $q_1\,,q_2\,,q_3$ be their numbers of arcs. Then $\xi=\pm(\xi_1-\xi_2)$ and $\zeta=\pm(\xi_3-\xi_2)$, so $\xi\times\zeta=\pm\bigl(\xi_1\times\xi_2+\xi_2\times\xi_3+\xi_3\times\xi_1\bigr)$, and by \eqref{eq:TS}, $\sum_{e\in S}(\xi\times\zeta)_e=\pm\bigl(\det(w_1\,,w_2)+\det(w_2\,,w_3)+\det(w_3\,,w_1)\bigr)$ with $w_i:=T_S\xi_i\in[-q_i\,,q_i]^2$ for $i=1,2,3$. This expression is affine in each $w_i$ separately, so its extreme values over the three boxes are attained when each $w_i$ is a corner, $w_i=q_i\varepsilon_i$ with $\varepsilon_i\in\{-1,1\}^2$. Among three elements of $\{-1,1\}^2$, two coincide or two are opposite, so one of the three determinants vanishes, while each of the others has absolute value at most twice the product of the corresponding path lengths. Hence
\[
\Bigl|\textstyle\sum_{e\in S}(\xi\times\zeta)_e\Bigr|\leq2\max\{q_1q_2+q_2q_3\,,\ q_2q_3+q_3q_1\,,\ q_3q_1+q_1q_2\}\leq\frac{(q_1+q_2+q_3)^2}{2}\leq\frac{(d-1)^2}{2},
\]
using $2q_2(q_1+q_3)\leq\bigl((q_1+q_3)+q_2\bigr)^2/2$ and its analogues. In both cases $W(\xi\times\zeta)=\max_S|\sum_{e\in S}(\xi\times\zeta)_e|\leq(d-1)^2/2$, which proves the theorem.
\end{proof}

\begin{corollary}\label{cor:constants}
$1/3\leq\liminf_{d\to\infty}H_3(d)/d^2\leq\limsup_{d\to\infty}H_3(d)/d^2\leq1/2$.
\end{corollary}

\begin{proof}
Combine \cref{cor:grid} and \cref{thm:upper3}.
\end{proof}

For $\ell=4$, \cref{thm:posetheight} together with $\delta_4=3$ (\cref{lem:delta4}) makes the leading constant explicit.

\begin{corollary}\label{cor:constants4}
Let $P$ be a poset on $[d]$ with $\ell(P)=4$. Every primitive facet normal $\alpha$ of $M(P)$ satisfies $\lVert\alpha\rVert_\infty\leq3\binom{d-1}{3}=(d-1)(d-2)(d-3)/2$. Consequently, $1/27\leq\liminf_{d\to\infty}H_4(d)/d^3\leq\limsup_{d\to\infty}H_4(d)/d^3\leq1/2$. \end{corollary}

\begin{proof}
The bound is \cref{thm:posetheight} with $\delta_4=3$ (\cref{lem:delta4}); it is $3/8$ of the bound $4\,\delta_3\binom{d-1}{3}=8\binom{d-1}{3}$ obtained from \cref{lem:posetbounds,lem:window}. The lower bound on the $\liminf$ is \cref{cor:growth} with $\ell=4$, which uses \cref{thm:product} with three nearly equal factors.
\end{proof}

\begin{remark}[Open problems]\label[remark]{rem:open}
For fixed $\ell\geq2$, $H_\ell(d)$, equivalently the largest $x$-coefficient of an inequality describing a vertical facet of $\D$ over posets with $\ell(P)=\ell$ (\cref{cor:bijection}), is $\Theta(d^{\ell-1})$ by \cref{cor:growth}; the constants lie in $[1/3,1/2]$ for $\ell=3$ and in $[1/27,1/2]$ for $\ell=4$ (\cref{cor:constants,cor:constants4}). In general, the constants of \cref{cor:growth} differ by the factor $\delta_\ell(\ell-1)^{\ell-1}/(\ell-1)!$, which is $1,4,27/2,160/3$ for $\ell=2,3,4,5$. The exact values of $H_\ell(d)$ for $\ell\geq3$ and the sharp constants remain open. Complete enumeration gives $H_3(5)=3$ and $H_3(6)=5$ (\cref{rem:sixteen}).

For $\ell\geq3$ and small $d$, the bounds of \cref{thm:posetheight,thm:upper3} are weaker than bounds that ignore the poset. Every facet normal of $M(P)$ is a facet normal of a Minkowski sum of staircase simplices, so its primitive representative has entries at most $\delta_{d-1}$ (\cref{rem:posetheight}), and at most the largest coefficient $\operatorname{coeff}(d)$ of a facet-describing inequality, with coprime integer coefficients, of a full-dimensional $0/1$-polytope in $\mathbb{R}^d$ \SF{\SFcoeff}; complete enumeration gives $\operatorname{coeff}(8)=18$ \cite[Section~5]{Ziegler01}. Since $\delta_6=9$ (see \cite{OEISA003432}), $H_\ell(7)\leq9$ and $H_\ell(8)\leq18$ for every $\ell$, whereas \cref{thm:upper3} gives only $H_3(7)\leq18$ and $H_3(8)\leq24$, and \cref{cor:constants4} gives $H_4(7)\leq60$. With \cref{cor:grid,prop:lowerboundfour}, this gives
\[
H_3(7)\in\{8,9\},\qquad H_4(7)=9,\qquad H_3(8)\in[11,18],\qquad H_4(8)\in[13,18].
\]
Thus the first unknown value is $H_3(7)$, and because its entries are at most $9$, it is decided by a finite search. Searches for $\ell=3$ and $d=7,8,9$ (over vectors with bounded entries for $d=7,8$, along facet adjacencies for $d=9$, each candidate tested by \cref{lem:dimofFMP}) found no entries above the lower bound $8,11,16$ of \cref{cor:grid}; we do not know whether that bound is sharp in general, and for $\ell=4$ and $7\leq d\leq9$, \cref{prop:lowerboundfour} beats the constructions of this section. A natural lead for $\ell=3$ is that the boxes in \eqref{eq:TS} are loose: $T_{\{1,2,3\}}$ maps $\{0,1\}^E$ onto a hexagon that omits the corners $(1,-1)$ and $(-1,1)$ of $[-1,1]^2$, so the estimate $2q_1q_2$ in the proof of \cref{thm:upper3} is not attained by label vectors.
\end{remark}

\medskip
\noindent {\bf Declaration of the use of AI in the manuscript preparation.}
During the preparation of this work, the authors used Claude to assist in developing and checking
the work. The authors take full responsibility for the content of the article.

\newpage

\bibliographystyle{siamplain}
\bibliography{orderpolytope}

\end{document}